\documentstyle{amsppt}
\magnification1200 
\NoBlackBoxes
\def\phi{\varphi}
\def\RE{\text{\rm Re }}

\def\lap{{\Cal L}}
\def\hbar{\overline{h}}

\def\Lam{\Lambda}
\def\tf{{\tilde f}}
\def\hchi{{\hat \chi}}
\def\hsigma{{\hat \sigma}}

\pagewidth{6.5 true in}
\pageheight{9 true in}

\topmatter
\title
Decay of mean-values of multiplicative functions 
\endtitle
\author
Andrew Granville and K. Soundararajan
\endauthor
\rightheadtext{Decay of mean-values of Multiplicative Functions}
\thanks{The first author is a Presidential Faculty Fellow.  He is also
supported, in part, by the National Science Foundation.  The second 
author is 
supported by the American Institute of Mathematics (AIM), 
and, in part, by the National Science Foundation (DMS 97-29992).}\endthanks

\address
Department of Mathematics, University of Georgia, Athens, GA , USA
\endaddress
\email 
andrew\@sophie.math.uga.edu
\endemail
\address
School of Mathematics,  Institute for Advanced Study,
Princeton, NJ 08540, USA
\endaddress 
\email 
ksound\@math.ias.edu
\endemail

\endtopmatter

\head 1. Introduction \endhead

\noindent Given a multiplicative function $f$ with $|f(n)|\le 1$ 
for all $n$, we are concerned with obtaining {\sl explicit} upper bounds on the 
mean-value $\frac{1}{x} |\sum_{n\le x} f(n)|$.  Ideally, one 
would like to give a bound for this mean-value which depends only on a 
knowledge of $f(p)$ for primes $p$.  To illustrate what we mean, we 
recall a pioneering result of E. Wirsing [16].  Throughout, we put 
$$
\Theta(f,x):= \prod_{p\le x} \Big(1+ \frac{f(p)}{p} +\frac{f(p^2)}{p^2} +\ldots
\Big) \Big(1-\frac 1p\Big).
$$
A. Wintner [15] showed by a simple convolution argument that if $\sum_p 
\frac{|1-f(p)|}{p}$ coverges then
$$
\lim_{x\to \infty} \frac{1}{x} \sum_{n\le x} f(n) = \Theta(f,\infty). \tag{1.1}
$$
If we restrict ourselves to real-valued multiplicative 
functions, 
then Wirsing showed that if $\sum_p \frac{1-f(p)}{p}$ diverges then
the limit in (1.1) exists, and equals $0 =\Theta (f,\infty)$.  Wirsing's result
settled an old conjecture of P. Erd{\H o}s and Wintner that 
every multiplicative function $f$ with $-1\le f(n)\le 1$ had a mean-value.

The situation for complex valued multiplicative functions is more delicate. 
For example, the function $f(n) = n^{i\alpha}$ ($0 \neq \alpha \in {\Bbb R}$)
does not have a mean-value because $\frac{1}{x} \sum_{n\le x} n^{i\alpha} 
\sim \frac{x^{i\alpha}}{1+i\alpha}$.  Note that here 
$\sum_p (1-\RE p^{i\alpha})/p$ diverges but $x^{-1} \sum_{n\le x} 
n^{i\alpha}$ does not tend to $0$.  G. Hal{\' a}sz [5, 6] realized that 
for complex valued multiplicative functions, the analogue of Wirsing's  
result requires the divergence of $\sum_p (1-\RE f(p)p^{-i\alpha})/p$ 
for all real numbers $\alpha$.  If this holds, then he showed that 
$\frac{1}{x}\sum_{n\le x} f(n) \to 0$, and he quantified how rapidly the 
limit is attained.  

\proclaim {Theorem (Hal{\'a}sz)} Let $f$ be a multiplicative function with 
$|f(n)|\le 1$ for all $n$, and set 
$$
M(x,T) = \min_{|y| \le 2T} \sum_{p\le x} 
\frac{1-\RE f(p) p^{-iy}}{p}. \tag{1.2} 
$$
Then 
$$
\frac{1}{x} \Big| \sum_{n\le x} f(n) \Big| \ll \exp\Big(-\frac{M(x,\frac 12 \log x)}{16}
\Big).
$$
\endproclaim 

Hal{\' a}sz comments that the factor $1/16$ may be replaced by the 
optimal constant $1$.  Our first Theorem provides such a refinement of 
Hal{\' a}sz' result:  Hal{\' a}sz' statement is a little 
inaccurate, the extra factor $M$ in our result below is necessary.

\proclaim{Theorem 1}  Let $f$ be a multiplicative function 
with $|f(n)| \le 1$ for all $n$.  Let $x\ge 3$, and let $T\ge 1$ be 
real numbers.  Put for any complex number $s$ with Re$(s)>0$,
$$
F(s)= \prod_{p\le x} \Big(1+\frac{f(p)}{p^s} + \frac{f(p^2)}{p^{2s}} +\ldots
\Big),
$$
and let 
$$
L = L(x,T) = \frac{1}{\log x} \Big( \max_{|y|\le 2T} |F(1+iy)|\Big). \tag{1.3}
$$
Then 
$$
\frac{1}{x}\Big|\sum_{n\le x} f(n) \Big| \le L \Big( \log \frac{e^{\gamma}}{L}
+ \frac{12}{7}\Big) + O\Big(\frac{1}{T}+\frac{\log \log x}{\log x}\Big).
$$
\endproclaim

\proclaim{Corollary 1}  Let $f$, $x$, and $T$ be as in Theorem 1, and let 
$M=M(x,T)$ as in (1.2). 
If $f$ is completely multiplicative then 
$$
\frac{1}{x} \Big| \sum_{n\le x} f(n)\Big|
\le \Big(M+\frac{12}{7}\Big) e^{\gamma -M} + O\Big( \frac{1}{T} + 
\frac{\log \log x}{\log x}\Big).
$$
If $f$ is multiplicative then 
$$
\frac{1}{x} \Big| \sum_{n\le x} f(n)\Big|
\le \prod_{p}\Big(1+ \frac{2}{p(p-1)}\Big) 
\Big(M+\frac{4}{7}\Big) e^{\gamma -M} + O\Big( \frac{1}{T} + 
\frac{\log \log x}{\log x}\Big).
$$
\endproclaim

As we will discuss after Theorem 5, Corollary 1 (and so Theorem 1) 
is essentially ``best possible'' (up to a factor $10$)
in that for any given $m_0$, 
we can construct $f$ and $x$ so that $M=M(x,\infty)>m_0$ and 
$|\sum_{n\le x} f(n)| \ge (M+12/7) e^{\gamma -M}/10$.

The maximum in (1.3) and the minimum in (1.2)) are a 
little unwieldly to compute, and it would 
be desirable to get similar 
decay estimates in terms of $|F(1)|$ (or, equivalently 
$\sum_{p\le x} (1-\RE f(p))/p$).  
In light of Hal\' asz's work (and particularly the example 
$f(n)=n^{i\alpha}$) this is possible only if we have some additional 
information on $f(n)$, such as knowing that all $f(p)\in D $ for some 
closed convex subset $D$ of the unit disc ${\Bbb U}$.  Such variants have been 
considered by Hal{\' a}sz [5,6], R. Hall and G. Tenenbaum [10], and Hall [9].  
The result of Hall is the most general, and to describe it we 
require some information on the geometry of $D$.  We collect this 
in Lemma 1.1 below, which is mostly contained in Hall's work.

\proclaim{Lemma 1.1}  Let $D$ be a closed, convex subset of ${\Bbb U}$ with $1
\in D$.  For $\alpha \in [0,1]$ define 
$$
\hbar(\alpha) = \frac{1}{2\pi} \int_0^{2\pi} \max_{\delta \in D} 
\RE (1-\delta)(\alpha -  e^{-i\theta}) d\theta. 
\tag{1.4}
$$
Define $\kappa = \kappa(D)$ to be the largest value of $\alpha\in [0,1]$ 
such that $\hbar(\alpha)\le 1$.  Lastly, put $\nu=\max_{\delta \in D} (1-\RE 
\delta)$.  Then, $\hbar$ is a continuous, increasing, convex function of 
$\alpha$, and $2\pi\hbar(0) = \lambda(D)$, the perimeter length of 
the boundary of $D$ (so $\kappa$ exists).  For $\kappa$ we have the lower bound 
$$
\kappa \ge \min\Big( 1, \frac{1-\hbar(0)}{\hbar(1)-\hbar(0)}\Big) 
\ge \min\Big(1, \frac{1}{\nu} \Big(1-\frac{\lambda(D)}{2\pi}\Big)\Big).
$$
Finally, $\kappa \nu \le 1$ for all $D$, and equality holds here 
if and only if $D = [0,1]$.  
\endproclaim

\remark{Remark} Hall also showed that if $0\in D$ then $\kappa=0$ only
when $D={\Bbb U}$, and  $\kappa=1$ only when $D=[0,1]$.
\endremark

\proclaim{Theorem (Hall)}  Retain the notations of Lemma 1.1, and 
let $f$ be a multiplicative function with $|f(n)|\le 1$, and 
$f(p)\in D$ for all primes $p$.  Then 
$$
\frac{1}{x }\Big| \sum_{n\le x} f(n)\Big| 
\ll  \exp\biggl(-\kappa(D) \sum_{p\le x} \frac{1-\text{Re }f(p)}{p}\biggr).
\tag{1.5} 
$$
\endproclaim

Hall states this result under the additional constraint that $0\in D$, 
but this is not necessary.  Hall also observed that the constant $\kappa(D)$ 
in (1.5) is optimal for {\sl every} $D$: it cannot be replaced by any larger value.  For 
completely multiplicative functions, we have obtained the following explicit 
version of Hall's theorem.

\proclaim{Theorem 2}  Retain the notations of Lemma 1.1, and further 
define 
$$
C(D) = -\kappa \nu \gamma + \min_{\epsilon = \pm 1} \int_0^{2\pi} 
\frac{\min(0,1-\kappa-\max_{\delta\in D} 
\RE \delta(e^{\epsilon ix}-\kappa))}{x} 
dx. \tag{1.6}
$$
Let $f$ be a multiplicative function with $f(p)\in D$ for all 
primes $p$, and put $y= \exp((\log x)^{\frac 23})$.  
If $\kappa \nu <1$ then 
$$ 
\align
\frac{1}{x} \Big|\sum_{n\le x} f(n)\Big| 
&\le  |\Theta(f,y)| \Big(\frac{2-\kappa \nu}{1-\kappa \nu}\Big)
 \exp\Big(-\kappa \sum_{y<p\le x} 
\frac{1-\RE f(p)}{p}  -C(D)+ \gamma(1-\kappa \nu)\Big) \\
&+ 
O\Big( \frac{1}{(\log x)^{\frac 13}} \exp\Big( \Big(2\log\log x \sum_{p\le x} 
\frac{1-\RE f(p)}{p} \Big)^{\frac 12}\Big)\Big).
\\
\endalign
$$
If $\kappa \nu =1$ (so that $D=[0,1]$) then 
$$
\frac{1}{x} \sum_{n\le x} f(n) \le e^{\gamma} |\Theta(f,x)|  
+ O\Big(\frac{1}{\log x}\Big).
$$
\endproclaim 

A version of the second statement in Theorem 2 was
first proved by Hall [9]. Theorem 2 
is essentially ``best possible'' (up to the constant of multiplication),
for {\sl every} such $D$, as noted in [9] and [10].

The first statement of Theorem 2 gives an explicit quantitative version of 
Hall's theorem, so long as $\sum_{p\le x} (1-\RE f(p))/p \ll \log \log x$.  
When $\sum_{p\le x} (1-\RE f(p))/p \gg \log \log x$, then the bound of 
Theorem 2 is no longer useful; and one should revert back to Hall's estimate.
In the case that $\sum_{p\le x} (1-\RE f(p))/p \gg \log \log x$, then 
Hall's theorem shows that $\frac{1}{x} \sum_{n\le x} f(n) = o(1)$.  
So we have the following Corollary to Theorem 2: 

\proclaim{Corollary 2}  Retain the notations of Theorem 2.  Let $f$ be 
a multiplicative function with $f(p)\in D$ for all $p$.  If $\kappa \nu 
<1$ then 
$$
\frac{1}{x}\Big| \sum_{n\le x} f(n) \Big| \le |\Theta(f,y)| 
\Big(\frac{2-\kappa \nu}{1-\kappa \nu}\Big)
 \exp\Big(-\kappa \sum_{y<p\le x} 
\frac{1-\RE f(p)}{p}  -C(D)+ \gamma(1-\kappa \nu)\Big) + o(1).
$$
\endproclaim

If the maximum in (1.3) (or, the minimum in (1.2)) occurs for $y=y_0$ then 
$f(n)$ looks roughly like $n^{iy_0}$, so that the mean-value of $f(n)$ 
should be of size $|\frac{x^{iy_0}}{1+iy_0}| \asymp \frac{1}{1+|y_0|}$.  
Our next result confirms this expectation.  

\proclaim{Theorem 3}  Let $f$ be a multiplicative function 
with $|f(n)|\le 1$ for all $n$.  Take $T=\log x$ in Theorem 1, and suppose 
the maximum in (1.3) is attained at $y=y_0$.  Then 
$$
\frac{1}{x} \Big| \sum_{n\le x} f(n) \Big| \ll \frac{1}{1+|y_0|} 
+ \frac{(\log \log x)^{1+2(1-\frac{2}{\pi})}}{(\log x)^{1-\frac{2}{\pi}}}.
$$
\endproclaim

Evidently this is ``best possible'', since taking $f(n)=n^{iy_0}$
gives the right side of the equation.

Lastly, we give an application of our ideas to the variation of averages of 
multiplicative functions.  Ideally, one would like to say that 
$$
\frac{1}{x} \sum_{n\le x} f(n) - \frac{w}{x} \sum_{n\le x/w} f(n) 
\ll \Big( \frac{\log 2w}{\log x}\Big)^\beta, 
\tag{1.7}
$$ 
for all $1\le w\le x$, with as large an exponent $\beta$ as possible 
($\beta=1$ would be optimal).  This would show that averages of multiplicative
functions vary slowly.  Unfortunately, (1.7) is not true in general, as 
the example $f(n)=n^{i\alpha}$ reveals.  However, P.D.T.A. Elliott [2] realized
that the absolute value of averages of multiplicative functions  
always varies slowly.  He showed that 
$$
\frac{1}{x} \Big| \sum_{n\le x} f(n) \Big| - \frac{w}{x} \Big| \sum_{n\le x/w}
f(n) \Big| \ll \Big( \frac{\log 2w}{\log x}\Big)^{\frac{1}{19}}, 
$$ 
for all multiplicative functions $f$ with $|f(n)|\le 1$, and all $1\le w\le x$.
One application of such an estimate, 
as Hildebrand [11] observed, is to (slightly)  extending
the range of validity of Burgess' character sum estimate.  
By applying Theorem 1, and the ideas underlying it, we have obtained 
the following improvement on Elliott's result.  
We remark that $1-\frac{2}{\pi} = 0.36338\ldots$, and $2-\sqrt{3} 
= 0.267949\ldots$.

\proclaim{Theorem 4} Let $f$, $x$, and $F$ be as in Theorem 1.  
Take $T=\log x$, and suppose that the maximum in (1.3) occurs 
at $y_0$.  
Then for $1\le w \le x/10$, we have
$$
\!\Big|
\frac{1}{x} \sum_{n\le x} f(n) n^{-iy_0} - \frac{w}{x} \sum_{n\le x/w} f(n) 
n^{-iy_0} \Big| 
\ll\! \Big(\frac{\log 2w}{\log x}\Big)^{1-\frac{2}{\pi}}\log \Big(\frac{\log x}
{\log 2w}\Big) + \frac{(\log \log x)^{1+2(1-\frac {2}{\pi})}}
{(\log x)^{1-\frac{2}{\pi}}}. 
$$
\endproclaim

\proclaim{Corollary 3} Let $f$ be a multiplicative function with $|f(n)|\le 
1$ for all $n$.  Then for $1\le w\le x/10$, we have
$$
\frac 1x \biggl| \sum_{n \le x} f(n) \biggr| - \frac wx 
\biggl| \sum_{n\le x/w}  f(n) \biggr| \ll \biggl(\frac{\log 2w}{\log x}
\biggr)^{1-\frac{2}{\pi}}\log \left( \frac{\log x}{\log 2w} \right)
 + \frac{\log \log x}{(\log x)^{2-\sqrt{3}}}.
$$
\endproclaim

Our proofs of Theorems 1, 3, and 4 are based on the following key Proposition
(and its variant Proposition 3.3 below),
which we establish by a variation of Hal{\' a}sz' method.  
Proposition 1 below is a variant of Montgomery's lemma (see [12], and also 
Montgomery and R.C. Vaughan [14]) which is 
one of the main ingredients in the proof of Hall's theorem.

\proclaim{Proposition 1}  Let $f$, $x$, $T$, and 
$F$ be as in Theorem 1.  
Then 
$$
\frac 1x \Big| \sum_{n\le x} f(n) \Big| \le \frac{2}{\log x} 
\int_0^{1} \Big(\frac{1 -x^{-2\alpha}}{2\alpha} \Big)
\Big(\max_{|y|\le T} |F(1+\alpha+iy)| \Big) d\alpha
+ O\Big(\frac {1}{T} + \frac{\log \log x}{\log x}\Big).
$$
\endproclaim

To prove Theorem 2, we adopt a different strategy, turning to integral 
equations. Let $\chi:[0,\infty) \to {\Bbb U}$ be a measurable function, with $\chi(t)
=1$ for $t\le 1$.  We let $\sigma(u)$ denote the solution to 
$$
\align
& \quad u\sigma(u) = (\sigma*\chi)(u) = \int_0^u \sigma(t) \chi(u-t) dt, 
\tag{1.8}\\
& \text{with initial condition } \sigma(u) = 1 \text{ for } 0\le u\le 1.
 \\
\endalign
$$
We showed in [4] that (1.8) has a unique solution, and this solution is 
continuous.  Further let $I_0(u;\chi)= 1$, and for $n\ge 1$ define 
$$
I_n(u;\chi) = \int\Sb t_1,\ldots,t_n\\ t_1+\ldots +t_n \le u\endSb 
\frac{1-\chi(t_1)}{t_1} \frac{1-\chi(t_2)}{t_2} \ldots \frac{1-\chi(t_n)}{t_n}
dt_1 \cdots dt_n. \tag{1.9a}
$$
Then we showed that 
$$
\sigma(u) = \sum_{n=0}^{\infty} \frac{(-1)^n}{n!} I_n(u;\chi). \tag{1.9b}
$$
The relevance of the class of integral equations (1.8) to the study of 
multiplicative functions was already observed by Wirsing [16].  We illustrate 
this connection by means of the following Proposition, proved in [4] 
(Proposition 1 there).

\proclaim{Proposition 2} Let $f$ be a multiplicative function with $|f(n)|
\le 1$ for all $n$ and $f(n)=1$ for 
$n\le y$.  Let $\vartheta(x) =\sum_{p\le x} \log p$ and define 
$$
\chi(u) = \chi_f(u) = \frac{1}{\vartheta(y^u) } \sum_{p\le y^u} f(p)\log p. 
$$
Then $\chi(t)$ is a measurable function taking values in the unit disc
and with $\chi(t)=1$ for $t\le 1$. Let 
$\sigma(u)$ be the corresponding unique solution to (1.8).
Then 
$$
\frac{1}{y^u} \sum_{n\le y^u} f(n) =\sigma(u) +O\biggl(\frac{u}{\log y}\biggr).
$$
\endproclaim

Proposition 2 allows us to handle mean-values of multiplicative functions
which are known to be $1$ on the small primes.  We borrow another result 
from [4] (see Proposition 4.5 there) which allows us to remove the 
impact of the small primes.

\proclaim{Proposition 3} Let $f$ be a  multiplicative function 
with $|f(n)|\leq 1$ for all $n$.  For any $2\le y\le x$, let $g$ be the
completely multiplicative function with $g(p)=1$ if $p\leq y$,
and $g(p)=f(p)$ otherwise. Then
$$
\frac 1x \sum_{n\le x} f(n) = \Theta (f,y)\ \frac 1x \sum_{m\le x} 
g(m) + O\Big( \frac{\log y}{\log x} \exp\Big(\sum_{p\le x} \frac{|1-f(p)|}{p}
\Big)\Big).
$$
\endproclaim

We prove Theorem 2 by establishing a decay estimate, Theorem 5,
 for solutions of (1.8)
when $\chi(t)$ is constrained to lie in $D$ for all $t$.  Then 
using Propositions 2 and 3 we unwind this result to deduce Theorem 2.

\proclaim{Theorem 5} Let $\chi:[0,\infty) \to D$ 
be a measurable function with $\chi(t)=1$ 
for $t\le 1$, and let $\sigma$ denote the corresponding solution to (1.8). 
Retain the notations of Lemma 1.1 and (1.5), and put 
$$
M_0 = M_{0}(u;\chi) = \int_0^u \frac{1-\RE \chi(v)}{v} dv.
$$
Then, if $\kappa\nu <1$,
$$
|\sigma(u)| \le \Big(\frac{2-\kappa \nu}{1-\kappa \nu}\Big)
 \exp\left(-\kappa M_0 -C(D)+ \gamma(1-\kappa \nu)\right)
 - 
\Big(\frac{\kappa \nu}{1-\kappa \nu}\Big)\exp\Big(-\frac{M_0}{\nu}-\frac{C(D)}
{\kappa \nu}\Big).
$$
If $\kappa \nu=1$ (so that $D=[0,1]$) then $|\sigma(u)|\le e^{\gamma -M_0}$.
\endproclaim

When studying mean values of multiplicative functions we have seen how the
example $f(n)=n^{i\alpha}$ led Hal\' asz to consider convex regions $D$
that are not dense on the unit circle. Given that we now have $\chi(t)=1$
for $0\leq t\leq 1$, it is perhaps unclear whether such restrictions
are necessary when considering (1.8). In fact they are, and in
section 10a we shall see that if $\chi(t)=e^{i\alpha t}$ for all $t>1$
then $\lim \sup |\sigma(u)|\gg_\alpha 1$.

By Proposition 2, we know that statements about 
multiplicative functions, can be interpreted to give information on solutions 
to (1.8).  For example, the remark after the statement of Theorem 2 translates
to saying that Theorem 5 is ``best possible'' for {\sl every} $D$, 
up to the constant of multiplication, via [9] and [10].
Moreover we can state
integral equations versions of Corollary 1 and Theorem 4.

\proclaim{Corollary 1$'$}  If $\chi$ and $\sigma$ are as in Theorem 5 then
$|\sigma(u)| \le (M+12/7) e^{\gamma -M}$ where
$$
M=M(u):= \min_{y\in {\Bbb R}} \int_0^u \frac{1-\RE \chi(v)e^{-ivy}}{v} dv .
$$
\endproclaim

In fact this is ``best possible'', up to a factor $10$,
 in the sense that for any given $m_0$
we can find $\chi$ and $\sigma$ as in Theorem 5 with $M>m_0$ and
$|\sigma(u)| \geq (M+12/7) e^{\gamma -M}/10$; see section 10b for our construction.
This implies the same of Corollary 1 and hence of Theorem 1.

The analogue of Theorem 4 shows that  $|\sigma(u)|$ 
obeys a strong Lipschitz-type estimate.

\proclaim{Theorem 4$'$}  Let $\chi:[0,\infty) \to {\Bbb U}$ be a measurable function
with $\chi(t)=1$ for $t\le 1$, and let $\sigma$ denote the 
corresponding solution to (1.8).  Then for all $1\le v\le u$, 
$$
\Big||\sigma(u)| -|\sigma(v)|\Big| \ll 
\Big(\frac{u-v}{u}\Big)^{1-\frac{2}{\pi}} \log \frac{u}{u-v}.
$$
\endproclaim

We illustrate Theorem 5, and thus Theorem 2,
by working out several examples. In each of our examples we will have 
$D=\overline{D}$, which allows us to restate Theorem 5 as
$|\sigma(u)|\leq c' e^{-\kappa M_0}<c e^{-\kappa M_0}$ where 
$$
c' := c \exp \left( -2\pi \int_0^{\pi} 
\frac{\min(0,1-\kappa-\max_{\delta\in D} \RE \delta(e^{i\theta }-\kappa))}{\theta (2\pi -\theta)} 
d\theta \right) \  < \ c:=\Big(\frac{2-\kappa \nu}{1-\kappa \nu}\Big)
e^\gamma .
$$

\noindent {\bf Example 1}.   {\sl $D$ is the convex hull of the $m$th roots of unity}.
For $m=2$ we have  $D=[-1,1]$, $\nu=2$,  $\kappa=0.32867416320\ldots$ and $c'=6.701842225\dots<c=6.978982\dots$.
For larger $m$ we can determine a formula for $\hbar(\alpha)$; for example,
for odd $m\geq 3$, define $\delta_j=\theta_j-\pi(2j-1)/m$ where 
$\sin \theta_j/(\cos \theta_j - \alpha) = \tan(\pi(2j-1)/m)$, for
$1\leq j\leq (m+1)/2$. Then
$$
\hbar(\alpha) = \alpha +  \frac{1}{\pi} \left( \sin \frac{\pi}{m}  \left(1 + 2\sum_{j=1}^{(m-1)/2} \cos \delta_j\right) -\alpha \left( \delta_1 + \sum_{j=1}^{(m-1)/2}
(\delta_{j+1}-\delta_j) \cos  \frac{2\pi j}{m} \right)
\right)
$$
An analogous formula holds for even $m$. 
We computed $\kappa$ and $c$ (not $c'$) for various $m$:

\smallskip
\centerline{\vbox{\offinterlineskip
\hrule
\halign{& \vrule \ # & \ \vrule \, \vrule\ # & \ \vrule \ #& \ \vrule \ #& \ \vrule \ #& \ \vrule \ #& \ \vrule \ #&\ \vrule \ #&\ \vrule \ # &\strut \  \vrule #\cr
$m$ \hfil &\hfil$3$ &\hfil$4$ &\hfil$5$ &\hfil$6$ &\hfil$7$ &\hfil$8$ &\hfil$9$ &\hfil$10$&\cr
\noalign{\hrule}\cr
\noalign{\hrule} 
$\kappa$\hfil  &\hfil.167216&\hfil.098589&\hfil.063565&\hfil.044673&\hfil.032971&\hfil.025359 &\hfil.020086&\hfil.016305&\cr
\noalign{\hrule}\cr
$c$ \hfil &\hfil4.15845&\hfil3.99959&\hfil3.79356&\hfil3.73689&\hfil3.68124&\hfil3.65731&\hfil3.63435&\hfil3.62219&\cr
\noalign{\hrule}\cr}
\hrule}}
 
\smallskip

\centerline{\sl The $c$ and $\kappa$ values for $D$, the convex hull of the $m$th roots of unity.}
\medskip

\noindent One can show that, as $m\to 
\infty$, we have $\kappa = \pi^2/6m^2 + O(1/m^4)$ and $c=2e^\gamma+O(1/m^2)$.
Therefore, following the
proof of Theorem 2 of [4] we have that if $x$ is sufficiently
large and $p$ is a prime $\equiv 1 \pmod m$, then there are at
least $\{ \pi_m+o(1)\} x$ integers $\leq x$ which are $m$th power
residues $\pmod p$, where 
$\pi_m \geq \exp( -\exp ( \{ 3/\pi^4+o(1)\} m^4\log m))$.
(It is shown in [4] that $\pi_m \leq \exp ( -\{ 1+o(1)\} m \log m)$,
and  that $\pi_2=.1715\dots$, the only $m$ for which the best possible
 value has been determined).
\smallskip

\noindent {\bf Example 2}.   {\sl $D$ is the disc going through $1$
with radius $r\leq 1$}. Note that $\kappa=0$ if $r=1$. We have the
(relatively) simple formula,
$$
\hbar(\alpha) = r\left( \alpha + \frac{1}{\pi}  \int_{\theta_0}^\pi | e^{i\theta} - \alpha | d\theta \right) ,
$$
so that $\kappa=1$ if $r\leq \pi/(\pi+4) = .43990084\dots$.
For various radii $r$, we computed $\kappa$ and $c$:

\smallskip
\centerline{\vbox{\offinterlineskip
\hrule
\halign{& \vrule \ # & \ \vrule \, \vrule\ # & \ \vrule \ #& \ \vrule \ #& \ \vrule \ #& \ \vrule \ #& \ \vrule \ #&\ \vrule \ #&\ \vrule \ # &\strut \  \vrule #\cr
$r$ \hfil &\hfil$.4399..$ &\hfil$.45$ &\hfil$.5$ &\hfil$.6$ &\hfil$.7$ &\hfil$.8$ &\hfil$.9$ &\hfil$.95$&\cr
\noalign{\hrule}\cr
\noalign{\hrule} 
$\kappa$\hfil  &\hfil1&\hfil.968330&\hfil.822168&\hfil.580480&\hfil.390142&\hfil.236024&\hfil.108183&\hfil.051957&\cr
\noalign{\hrule}\cr
$c$ \hfil &\hfil16.5986&\hfil15.6413&\hfil11.7966&\hfil7.65099&\hfil5.70586&\hfil4.64287&\hfil3.99284&\hfil3.75723&\cr
\noalign{\hrule}\cr}
\hrule}}
 
\smallskip

\centerline{\sl The $c$ and $\kappa$ values for $D$, the disc of radius $r$, with center $1-r$.}
\medskip

\noindent One can show that, as $r$ gets close to $1$, that is $r=1-\delta$
where $\delta\to 0^+$, then  $\kappa = \delta + 3\delta^2/4 + O(\delta^3)$ and $c=2e^\gamma(1+\delta+O(\delta^2))$.
\smallskip

\noindent {\bf Example 3}.   {\sl $D$ is the sector of the circle bounded
by the lines from $1$ to $e^{\pm i\phi}$}. In other words,
$D$ is the convex hull of the point
$\{ 1\}$ together with the arc from $e^{i\phi}$ to $e^{-i\phi}$
on the unit circle.  
Select $\theta_0<\theta_1$ so that 
$\tan (\phi/2) = \sin \theta_0/(\cos \theta_0 - \alpha)$ and
$\tan \phi = \sin \theta_1/(\cos \theta_1 - \alpha)$, and thus,
with $I:=( \theta_0 + (\theta_1-\theta_0)\cos \phi)$, we have
$$
\hbar(\alpha) = \alpha + \frac{1}{\pi} \left(  
\sin \theta_0 + \sin(\theta_1-\phi) - \sin(\theta_0-\phi)
-\alpha I +
\int_{\theta_1}^\pi | e^{i\theta} - \alpha | d\theta \right) .
$$
Notice that if $\phi=\pi$ then $D=[-1,1]$ so, as above,
 $\kappa=\kappa^*:=.328674163\dots$ and $c=6.978982\dots$ We computed
the following values:
\smallskip
\centerline{\vbox{\offinterlineskip
\hrule
\halign{& \vrule \ # & \ \vrule \, \vrule\ # & \ \vrule \ #& \ \vrule \ #& \ \vrule \ #& \ \vrule \ #& \ \vrule \ #&\ \vrule \ #&\ \vrule \ # &\strut \  \vrule #\cr
$\phi$ \hfil &\hfil$\pi/4$ &\hfil$\pi/3$ &\hfil$\pi/2$ &\hfil$2\pi/3$ &\hfil$3\pi/4$ &\hfil$5\pi/6$ &\hfil$9\pi/10$ &\hfil$.99\pi$&\cr
\noalign{\hrule}\cr
\noalign{\hrule} 
$\kappa$\hfil  &\hfil.006293&\hfil.014597&\hfil.046181&\hfil.140280&\hfil.188459&\hfil.235961&\hfil.317918&\hfil.328674&\cr
\noalign{\hrule}\cr
$c$ \hfil &\hfil3.58485&\hfil3.61571&\hfil3.74339&\hfil4.01647&\hfil4.25671&\hfil 4.63956&\hfil5.15381&\hfil6.67192&\cr
\noalign{\hrule}\cr}
\hrule}}
 
\smallskip

\centerline{\sl The $c$ and $\kappa$ values for $D$, the cone with lines from $1$ to $e^{\pm i\phi}$.}
\medskip

\noindent One can show that as $\phi\to 0$ we have $\kappa\sim \phi^3/24\pi$.
Moreover if $\phi\to \pi$ then we have $\kappa^* - \kappa \sim
\eta (\pi-\phi)$, for some absolute constant $\eta$.

\head 2. Preliminaries \endhead 

\noindent We begin with the following lemma, weaker versions 
of which may be found in the works of Hal{\' a}sz [5], Halberstam and 
Richert [7], and Montgomery and Vaughan [14].   

\proclaim{Lemma 2.1} Let $f$ be a multiplicative function with 
$|f(n)| \le 1$ for all $n$.  Put $S(x) = \sum_{n\le x} f(n)$.  
Then for $x\ge 3$,
$$
|S(x)| \le \frac{x}{\log x} \int_2^x \frac{|S(y)|}{y^2}  dy 
+ O\Big(\frac{x}{\log x}\Big). \tag{2.1}
$$
Further, if $1\le w\le x$, then 
$$
\Big|\frac{S(x)}{x} - \frac{S(x/w)}{x/w}\Big| 
\le \frac{1}{\log x} \int_{2w}^x \Big| \frac{S(y)}{y} -\frac{S(y/w)}{y/w}\Big| 
\frac{dy}{y} + O\Big( \frac{\log 2w}{\log x}\Big). 
\tag{2.2}
$$
\endproclaim

\demo{Proof} First note that
$$
S(x)\log x   -\sum_{n \le x} f(n) \log n = 
\sum_{n \le x} 
f(n) \log \frac xn = O\biggl( 
\sum_{n \le x} \log \frac xn \biggr) = O(x). 
$$
Further
$$
\sum_{n\le x} f(n) \log n = \sum_{n\le x} f(n) \sum_{p^k| n} \log p 
=  \sum_{p^k \le x} \log p \sum_{m\le x/p^k} f(mp^k).
$$
Since 
$$
\sum_{m \le x/p^k} f(mp^k) = f(p^k) \sum_{m\le x/p^k} f(m) + O \biggl( 
\sum\Sb m \le x/p^k \\ p | m\endSb 1 \biggr) 
= f(p^k) S\Big(\frac{x}{p^k}\Big) + O\Big(\frac{x}{p^{k+1}}\Big),
$$ 
it follows that 
$$
S(x)\log x = \sum_{d\le x} f(d) \Lam(d)  S\Big(\frac xd\Big)
+ O(x).  \tag{2.3}
$$
Hence
$$
|S(x)|\log x 
\le \sum_{d \le x}\Lam(d)\Big|S\Big(\frac{x}{d}\Big)\Big| + O(x).
\tag{2.4}
$$

Writing $\psi(x) =\sum_{n\le x} \Lam(n)$, as usual, we see that 
$$
\sum_{d\le x} \Lam(d) \Big| S\Big(\frac xd\Big)\Big| 
= \sum_{d\le x} (\psi(d)-\psi(d-1)) \Big|S\Big(\frac xd\Big)\Big| 
= \sum_{d\le x} \psi(d) \Big( \Big|S\Big(\frac{x}{d}\Big)\Big| 
- \Big|S\Big(\frac{x}{d+1}\Big)\Big|\Big).
$$
We now use the prime number theorem in the form $\psi(d) 
= d+ O(d/(\log 2d)^2)$, together with the simple observation that 
$|S(x/d)| - |S(x/(d+1))| \le \sum_{x/(d+1) < n\le x/d} 1$.  
It follows that 
$$
\sum_{d\le x} \Lam(d) \Big| S\Big(\frac xd\Big)\Big| 
= \sum_{d\le x} d \Big( \Big|S\Big(\frac{x}{d}\Big)\Big| 
- \Big|S\Big(\frac{x}{d+1}\Big)\Big|\Big) + O\Big( \sum_{d\le x} 
\frac{d}{\log^2 (2d)} \sum_{x/(d+1) < n \le x/d} 1\Big).
$$
The main term above is plainly $\sum_{d\le x} |S(x/d)|$, and the 
remainder term is 
$$
\ll \sum_{d\le \sqrt{x}} \frac{d}{\log^2 (2d)} \frac{x}{d(d+1)} 
+ \frac{1}{\log^2 x} \sum_{\sqrt{x} \le d\le x} d \sum_{x/(d+1)<n \le x/d} 1
\ll x + \frac{1}{\log^2 x}\sum_{n\le \sqrt{x}} \frac{x}{n} \ll x.
$$
Combining these observations and (2.4), we have shown that 
$$
|S(x)| \log x \le \sum_{d\le x} \Big|S\Big(\frac{x}{d}\Big)\Big| + O(x).
$$
Now $|S(x/d)| = \int_d^{d+1} |S(x/t)| dt + O(\sum_{x/(d+1)< n \le x/d} 1)$,
and so the right side above is 
$$
\int_1^{x+1} \Big|S\Big(\frac{x}{t}\Big)\Big| dt + O(x).
$$
By changing variables $y=x/t$ this is 
$$
x \int_{x/(x+1)}^x \frac{|S(y)|}{y^2} dy +O(x) = x\int_2^x \frac{|S(y)|}{y^2}
dy 
+O(x),
$$
proving (2.1).

To show (2.2), we note by (2.3) that 
$$
\align
\log x \Big( \frac{S(x)}{x} -\frac{S(x/w)}{x/w} \Big) 
&= O(\log 2w) + \frac{1}{x} \sum_{d\le x} f(d)\Lam(d) S\Big(\frac{x}{d}\Big)
- \frac{w}{x} \sum_{d\le x/w} f(d) \Lam(d) S\Big(\frac{x}{wd}\Big)\\
&= O(\log 2w) 
+ \sum_{d\le x/w} f(d)\Lam(d) \Big(\frac{S(x/d)}{x} -\frac{S(x/wd)}{x/w}\Big).
\\
\endalign
$$
Hence
$$
\Big| \frac{S(x)}{x} - \frac{S(x/w)}{x/w}\Big| 
\le \frac{1}{\log x} 
\sum_{d\le x/w} \Lam(d) \Big| \frac{S(x/d)}{x} -\frac{S(x/wd)}{x/w}\Big|
+ O\Big(\frac{\log 2w}{\log x}\Big).
$$
We now mimic the partial summation argument used to deduce (2.1) from (2.4).  
This shows (2.2).

\enddemo

\proclaim{Lemma 2.2}  Let $a_n$ be a sequence of complex numbers such 
that $\sum_{n=1}^{\infty} \frac{|a_n|}{n} < \infty$.  Define 
$A(s) = \sum_{n=1}^{\infty} a_n n^{-s}$ which is absolutely convergent 
in Re$(s)\ge 1$.  For all real numbers $T \ge 1$, and all $0\le \alpha \le 1$ 
we have 
$$
\max_{|y|\le T} |A(1+\alpha+iy)| \le \max_{|y|\le 2T} |A(1+iy)| + O\Big(\frac 
\alpha T \sum_{n=1}^{\infty} \frac{|a_n|}{n}\Big), \tag{2.5}
$$
and for any $w\ge 1$,
$$
\max_{|y|\le T} |A(1+\alpha+iy)(1-w^{-\alpha-iy})| \le \max_{|y|\le 2T} 
|A(1+iy)(1-w^{-iy})| + O\Big( \frac \alpha T \sum_{n=1}^{\infty} \frac{|a_n|}{n}
\Big). 
\tag{2.6}
$$
\endproclaim
\demo{Proof}  We shall only prove (2.6); the proof of (2.5) is similar.  
Note that the Fourier transform 
of $k(z)= e^{-\alpha |z|}$ is ${\hat k}(\xi) = 
\int_{-\infty}^{\infty} e^{-\alpha |z| - i\xi z} dz = \frac{2\alpha}{\alpha^2 
+\xi^2}$ which is always non-negative.  The Fourier inversion formula gives
for any $z\ge 1$,
$$
z^{-\alpha} = k(\log z) = k(-\log z) = \frac{1}{2\pi} \int_{-\infty}^{\infty} 
{\hat k}(\xi) z^{-i\xi} d\xi = \frac{1}{\pi} \int_{-T}^{T} 
\frac{\alpha}{\alpha^2 +\xi^2} z^{-i\xi} d\xi + O\Big( \frac{\alpha}{T}\Big).
$$
Using this appropriately, we get that for all $n\ge 1$, and $0\le \alpha\le 1$,
$$
\frac{1}{n^{\alpha}} (1-w^{-\alpha-iy}) = \frac{1}{\pi} \int_{-T}^{T} 
\frac{\alpha}{\alpha^2+\xi^2} n^{-i\xi}(1-w^{-iy-i\xi}) d\xi + 
O\Big( \frac{\alpha}{T}\Big).
$$
Multiplying the above by $a_n/n^{1+iy}$, and summing over all $n$, we conclude 
that 
$$
A(1+\alpha+iy)(1-w^{-\alpha-iy}) = \frac{1}{\pi} \int_{-T}^{T} 
\frac{\alpha}{\alpha^2+\xi^2} A(1+iy+i\xi) (1-w^{-iy-i\xi}) d\xi 
+ O\Big(\frac{\alpha}{T} \sum_{n=1}^{\infty} \frac{|a_n|}{n}\Big).
$$
If $|y|\le T$ then $|y+\xi|\le |y|+|\xi| \le 2T$, and so 
we deduce that 
$$
\align
\max_{|y|\le T} |A(1+\alpha+iy)(1-w^{-\alpha-iy})| &\le 
\Big( \max_{|y|\le 2T} |A(1+iy)(1-w^{-iy})|\Big) \frac{1}{\pi} \int_{-T}^{T} 
\frac{\alpha}{\alpha^2+\xi^2}d\xi \\
&\hskip 1 in + O\Big( \frac \alpha T \sum_{n=1}^{\infty} 
\frac{|a_n|}{n}\Big),\\
\endalign
$$
and (2.6) follows since $\frac 1\pi \int_{-T}^T \frac{\alpha}{\alpha^2+\xi^2}
d\xi \le \frac{1}{\pi } \int_{-\infty}^{\infty} \frac{\alpha}{\alpha^2+\xi^2} 
d\xi = 1$.
\enddemo

Our next lemma was inspired by Lemma 2 of Montgomery and Vaughan [14],   
who consider (essentially) the quotient $|F(1+i(y+\beta))/F(1+iy)|$
rather than the product below.  

\proclaim{Lemma 2.3}  Let $f$,  $x$, and $F$ be as in 
Theorem 1.  Then for all real numbers $y$, and $1/\log x \le |\beta| 
\le \log x$, we have 
$$
|F(1+iy) F(1+i(y+\beta))| \ll (\log x)^{\frac{4}{\pi}} \max
\Big(\frac{1}{|\beta|}, (\log \log x)^2\Big)^{2(1-\frac 2{\pi})}.
$$
\endproclaim
\demo{Proof}  Clearly
$$
\align
|F(1+iy)F(1+i(y+\beta))| &\ll \exp\Big( \RE \sum_{p\le x} \frac{f(p)p^{-iy} 
+ f(p)p^{-i(y+\beta)}}{p}\Big)\\ 
&\ll \exp\Big( \sum_{p\le x} \frac{|1+p^{-i\beta}|}{p} \Big)=
\exp\Big(\sum_{p\le x}
\frac{2|\cos(\frac{|\beta|}{2} \log p)|}{p}\Big).
\tag{2.7}\\
\endalign
$$
By the prime number theorem and partial summation we have for $z\ge w\ge 2$ 
$$
\sum_{w\le p\le z} \frac{1}{p} = \int_w^z \frac{dt}{t\log t} 
+O(\exp(-c\sqrt{\log w})),
$$
for some constant $c>0$.  Choose $C=100/c^2$, and put $Y=
\max( \exp(C(\log \log x)^2), e^{\frac1{|\beta|}})$.  Put $\delta=1/\log^3 x$,
and divide the 
interval $[Y,x]$ into $\ll \log^4 x$ subintervals of the 
type $(z,z(1+\delta)]$ (with perhaps one shorter interval).  For each of these
subintervals we have
$$
\align
\sum_{z\le p \le z(1+\delta)} \frac{|\cos(\frac{|\beta|}{2}\log p)|}{p} 
&= ( |\cos(\tfrac{|\beta|}{2}\log z)| +O(\delta |\beta|)) 
\sum_{z\le p\le (1+\delta)z} \frac{1}{p} \\
&= ( |\cos(\tfrac{|\beta|}{2}\log z)| +O(\delta |\beta|)) \Big(  
\int_{z}^{z(1+\delta)}\frac{dt}{t\log t} + O\Big( \frac{1}{\log^{10} x}\Big)
\Big)
\\
&=\int_{z}^{z(1+\delta)} \frac{|\cos (\frac{|\beta|}{2} \log t)|}{t\log t} 
dt + O\Big( \frac{1}{\log^5 x}\Big),\\
\endalign
$$
where we used $|\beta| \le \log x$.
Using this for each of the $\ll \log^4 x$ such subintervals covering $[Y,x]$,
we conclude that 
$$
\sum_{Y\le p \le x} \frac{|\cos (\frac{|\beta|}{2}\log p)|}{p} 
= \int_{Y}^x \frac{|\cos (\frac{|\beta|}{2}\log t)|}{t\log t} dt 
+O\Big( \frac{1}{\log x}\Big) = \int_{\frac{|\beta|}{2} 
\log Y}^{\frac{|\beta|}{2}\log x} 
\frac{|\cos y|}{y} dy +O(1).
$$
Splitting the integral over $y$ above into intervals of length $2\pi$ 
(with maybe one shorter interval), and noting that $\frac{1}{2\pi} 
\int_0^{2\pi} |\cos\theta| d\theta =\frac{2}{\pi}$, we deduce that 
$$
\sum_{Y\le p \le x} \frac{|\cos (\frac{|\beta|}{2}\log p)|}{p} 
\le \frac{2}{\pi}
\log \frac{\log x}{\log Y} + O(1).
$$
Trivially, we also have 
$$
\sum_{p\le Y} \frac{|\cos (\frac{|\beta|}{2}\log p)|}{p} \le \sum_{p\le Y} 
\frac{1}{p} = \log \log Y +O(1).
$$
Combining the above two bounds,
 we get that 
$$
\align
\sum_{p\le x} \frac{|\cos (\frac{|\beta|}{2}\log p)|}{p} 
&\le \frac{2}{\pi} \log \log x + \Big(1-\frac{2}{\pi}\Big) \log \log Y +O(1) .
\\
\endalign
$$
The Lemma follows upon using this in (2.7), and recalling the definition 
of $Y$.
\enddemo

We conclude this section by offering a proof of Lemma 1.1.  

\demo{Proof of Lemma 1.1} For a fixed $\theta$, note that 
 $\max_{\delta \in D} \RE (1-\delta)(\alpha-e^{-i\theta})$ 
is an increasing function of $\alpha$.  Integrating, we see that 
$\hbar(\alpha)$ is an increasing function.  Clearly $\hbar$ is continuous, 
and we now show that it is convex: that is, given $0\le \alpha <\beta 
\le 1$, and $t\in [0,1]$, $\hbar(t\alpha+(1-t)\beta)\le t\hbar(\alpha)
+(1-t)\hbar(\beta)$.  Indeed, for a fixed $\theta$, we have
$$
\align
\max_{\delta \in D} \RE (1-\delta)(t(\alpha-e^{-i\theta}) + 
(1-t)(\beta-e^{-i\theta})) &\le t \max_{\delta \in D} \RE (1-\delta)(\alpha 
-e^{-i\theta}) \\
&+ (1-t) \max_{\delta \in D} \RE (1-\delta)(\beta -e^{-i\theta});
\\
\endalign
$$
so, integrating this, we get that $\hbar$ is convex.  

Note that $2\pi \hbar(0) = \int_0^{2\pi} \max_{\delta \in D} 
\RE (1-\delta)(-e^{-i\theta}) d\theta = \int_0^{2\pi} \max_{\delta \in D}
\RE \delta e^{-i\theta} d\theta$.  This last expression equals $\lambda(D)$,
the perimeter of $D$, a result known as Crofton's formula (see [1], page 65).

We now show the lower bounds for $\kappa$.  If $\kappa=1$ there is 
nothing to prove; and suppose $\kappa <1$ so that $\hbar(1)>1$.  
By convexity we see that 
$$
\hbar\Big(\frac{1-\hbar(0)}{\hbar(1)-\hbar(0)}\Big) 
\le \frac{1-\hbar(0)}{\hbar(1)-\hbar(0)}\hbar(1) + \Big(1- \frac{1-\hbar(0)}
{\hbar(1)-\hbar(0)}\Big) \hbar(0) = 1,
$$
and so it follows that $\kappa \ge \frac{1-\hbar(0)}{\hbar(1)-\hbar(0)}$.  
Clearly $\hbar(\alpha) \le \hbar(0) + \frac{1}{2\pi} \int_0^{2\pi} 
\max_{\delta \in D} \RE (1-\delta)\alpha d\theta = \hbar(0) +\alpha\nu$.  
Hence we see that $\frac{1-\hbar(0)}{\hbar(1)-\hbar(0)} \ge (1-\hbar(0))/\nu
= \frac{1}{\nu}(1-\frac{\lambda(D)}{2\pi})$.

Lastly it remains to show that $\kappa \nu \le 1$ with equality only when 
$D=\kappa \nu$.  By definition we have 
$\hbar(\alpha) \ge \max_{\delta \in D} 
\frac{1}{2\pi} \int_0^{2\pi} \RE (1-\delta)(\alpha-e^{-i\theta})d\theta
= \alpha \nu$.  It follows that $\kappa \nu \le 1$ always. Moreover if
   $\kappa \nu = 1$ then $\hbar(\kappa)=1$ and there exists $d\in D$ such
that the maximum of $\RE (1-\delta)(\kappa-e^{-i\theta})$ for $\delta\in D$,
occurs at $\delta=d$. Therefore if
$d+\eta \in D$ then $\RE \eta (\kappa-e^{-i\theta})\geq 0$ for all 
$\theta \in [0,2\pi)$. Therefore $\kappa=\nu=1$ else as $\theta$ runs through
$[0,2\pi)$, so does $\arg (\kappa-e^{-i\theta})$, which implies 
 $\RE \eta (\kappa-e^{-i\theta})< 0$ for some $\theta$.  Now 
$\arg (1-e^{-i\theta})$ runs through $(-\pi/2,\pi/2)$ so $\eta \in {\Bbb R}$
else $\RE \eta (\kappa-e^{-i\theta})< 0$ for some $\theta$. Thus
$D\subset {\Bbb R}$ and so $D=[0,1]$ since $\nu=1$.
\enddemo

\head 3. The key Proposition \endhead 

\subhead 3a.  The integral equations version \endsubhead 

\noindent Our tool in analysing (1.8) is the Laplace transform, which, for a 
measurable   function $f:\ [0,\infty) \to {\Bbb C}$ is given by
$$
{\Cal L}(f,s) = \int_0^{\infty} f(t) e^{-ts} dt
$$
where $s$ is some complex number.   If $f$ is integrable and 
grows sub-exponentially  
(that is, for every $\epsilon >0$, $|f(t)| 
\ll_{\epsilon} e^{\epsilon t}$ almost everywhere) then the 
Laplace transform  is well defined for all complex numbers 
$s$ with $\RE(s) >0$.  
Laplace transforms occupy a role in the study of differential equations 
analogous to Dirichlet series in multiplicative number theory.

Below, $\chi$ will be 
measurable with $\chi(t) =1$ for $t\le 1$ and $|\chi(t)| \le 1$ for all 
$t$, and $\sigma(u)$ will denote the corresponding solution to (1.8).
Observe that for any two `nice' functions $f$ and $g$, $\lap (f*g,s) =\lap 
(f,s)\lap(g,s)$.  From the definition of $\sigma$, it follows that 
$$
\lap(v\sigma(v),t+iy)= \lap(\sigma,t+iy) \lap(\chi,t+iy), \tag{3.1}
$$
where $t>0$ and $y$ are real numbers.  

Further, recalling from (1.9a,b)
that $\sigma(v)=\sum_{j=0}^{\infty} (-1)^j I_j(v;\chi)/j!$, we
have
$$
\align
\lap(\sigma,t+iy) &= \sum_{j=0}^{\infty} \frac{(-1)^j}{j!} \lap(I_j(v;\chi),
t+iy)= \frac{1}{t+iy} 
\sum_{j=0}^{\infty} \frac{(-1)^j}{j!} \biggl(\lap\biggl(\frac{1-\chi(v)}{v},
t+iy\biggr)\biggr)^{j} \\
&=\frac{1}{t+iy} \exp\biggl(-\lap\biggl(\frac{1-\chi(v)}{v},t+iy\biggr)\biggr).
\tag{3.2}\\
\endalign
$$
We now give our integral equations version of Proposition 1.

\proclaim{Proposition 3.1}  Let $u\ge 1$, and define
for $t>0$
$$
M(t) = \int_u^{\infty} \frac{e^{-tv}}{v} dv + \min_{y\in {\Bbb R}} 
\int_0^u \frac{1-\RE \chi(v)e^{-ivy}}{v} e^{-tv}dv.
$$
Then  
$$
|\sigma(u)| \le \frac{1}{u} \int_{0}^{\infty} 
\biggl(\frac{1-e^{-2tu}}{t}\biggr) 
\frac{\exp(-M(t))}{t}dt.
$$
\endproclaim

\noindent Since $M(t) \ge \max(0,-\log (tu)+O(1))$ we see that the 
integral in the Proposition converges.  

\demo{Proof}  Define $\hchi(v)=\chi(v)$ if $v\le u$, and $\hchi(v)=0$ if $v>u$.
Let $\hsigma$ denote the corresponding solution to (1.8).  Note that
$\hsigma(v)=\sigma(v)$ for $v\le u$.  Thus 
$$
\align
|\sigma(u)| &=|\hsigma(u)| \le \frac 1u \int_0^u |\hsigma(v)|dv 
= \frac{1}{u}\int_0^{u} 2v|\hsigma(v)| \int_0^{\infty} e^{-2t v}dt dv \\
&= \frac{1}{u}
\int_{0}^{\infty} \biggl(\int_0^{u} 2v|\hsigma(v)| e^{-2tv} dv 
\biggr) dt. \tag{3.3}\\
\endalign
$$

By Cauchy's inequality
$$
\align
\biggl(\int_{0}^{u} 2v |\hsigma(v)|e^{-2tv}dv \biggr)^2 &\le 
\biggl(4\int_{0}^{u} e^{-2tv} dv \biggr) 
\biggl(\int_0^{\infty} |v\hsigma(v)|^2 e^{-2tv} dv\biggr)\\
&= 2\frac {1-e^{-2tu}}t \int_0^{\infty} |v\hsigma(v)|^2 e^{-2tv} dv.\tag{3.4}\\
\endalign
$$
By Plancherel's formula (Fourier transform is an isometry on $L^2$)
$$
\int_0^{\infty} |v\hsigma(v)|^2 e^{-2tv} dv  = \frac{1}{2\pi}
\int_{-\infty}^{\infty} |\lap(v\hsigma(v),t+iy)|^2 dy 
$$
and, using (3.1), this is
$$
= \frac{1}{2\pi}\int_{-\infty}^{\infty} |\lap(\hsigma, t+iy)|^2 
|\lap(\hchi,t+iy)|^2 dy
\le \left( \max_{y\in {\Bbb R}} |\lap(\hsigma, t+iy)|^2\right)
 \frac{1}{2\pi}\int_{-\infty}^{\infty} |\lap(\hchi,t+iy)|^2dy.
$$
Applying Plancherel's formula again, we get
$$
\frac{1}{2\pi}
\int_{-\infty}^{\infty} |\lap(\hchi,t+iy)|^2 dy = \int_0^{\infty} |\hchi(v)|^2 
e^{-2tv} dv \le \int_0^{u} e^{-2tv} dv =\frac{1-e^{-2tu}}{2t}.
$$
Hence 
$$
\int_0^{\infty} |v\hsigma(v)|^2 e^{-2tv} dv \le \frac{1-e^{-2tu}}{2t} 
\max_{y\in {\Bbb R}} |\lap (\hsigma,t+iy)|^2. \tag{3.5}
$$

By (3.2), we have
$$
\lap(\hsigma,t+iy) 
=\frac{1}{t+iy} \exp\biggl(-\lap\biggl(\frac{1-\hchi(v)
e^{-ivy}}{v},t\biggr) + \lap\biggl(\frac{1-e^{-ivy}}{v},t\biggr)\biggr).
$$
Now, we have the identity 
$$
\RE \lap\biggl(\frac{1-e^{-ivy}}{v},t\biggr) = \log |1+iy/t|
$$
which is easily proved by differentiating both sides with respect to $y$.  
Using this we obtain 
$$
t|\lap(\hsigma, t+iy) | 
= \exp\biggl(-\RE \lap\biggl(\frac{1-\hchi(v)e^{-ivy}}{v}
,t\biggr)\biggr), \tag{3.6}
$$
from which it follows that 
$$
\max_{y\in {\Bbb R}} |\lap (\hsigma,t+iy)| = \frac{\exp(-M(t))}{t}.
$$
Inserting this in (3.5), and that into (3.4), and then (3.3), we 
obtain the Proposition.

\enddemo

\subhead 3b.  The multiplicative functions version: Proof of Proposition 1
\endsubhead

\noindent In this subsection, we prove Proposition 1.  We follow closely
 the ideas behind the proof of Proposition 3.1 above.  

Note that 
$$
\align
\int_2^x \frac{|S(y)|}{y^2}dy &= \int_2^x \frac{2\log y}{y^2} \Big| 
\sum_{n\le y} f(n)\Big| \Big( \int_0^{1} y^{-2\alpha} d\alpha +O(y^{-2})\Big)
dy 
\\
&= 
\int_2^x \frac{2}{y^2} \Big| \sum_{n\le y} f(n) \log n + O\Big(\sum_{n\le y} 
\log (y/n)\Big)\Big| \Big( \int_0^{1} y^{-2\alpha} d\alpha \Big)dy 
+O(1)\\
&= \int_0^{1} \Big(\int_2^x \frac{2}{y^{2+2\alpha}} \Big| \sum_{n\le y} 
f(n)\log n\Big| dy \Big)d\alpha + O( \log \log x). \tag{3.7}\\ 
\endalign
$$
By Cauchy's inequality 
$$
\align
\int_2^x \Big|\sum_{n\le y} f(n) \log n \Big| \frac{dy}{y^{2+2\alpha} } 
&\le \Big( \int_1^x \frac{dy}{y^{1+2\alpha}} \Big)^{\frac 12} 
\Big(\int_2^x \Big|\sum_{n\le y} f(n) \log n \Big|^2 \frac{dy}{y^{3+2\alpha}}
\Big)^{\frac{1}{2}}\\
&= \Big(\frac{1-x^{-2\alpha}}{2\alpha}\Big)^{\frac{1}{2}} 
\Big(\int_2^x \Big|\sum_{n\le y} f(n) \log n \Big|^2 \frac{dy}{y^{3+2\alpha}}
\Big)^{\frac{1}{2}}. \tag{3.8}\\
\endalign
$$

Now define the multiplicative function ${\tilde f}$ by ${\tilde f}(p^k) 
= f(p^k)$ for $p\le x$, and ${\tilde f}(p^k)=0$ for $p>x$, so that $F(s)=\sum_{n\geq 1}{\tilde f}(n)/n^s$.  Naturally 
$\tf (n) = f(n)$ for $n\le x$, and so 
$$
\int_2^x \Big|\sum_{n\le y} f(n) \log n \Big|^2 \frac{dy}{y^{3+2\alpha}}
\le \int_1^{\infty} \Big| \sum_{n\le y} \tf(n)\log n\Big|^2 \frac{dy}
{y^{3+2\alpha}},
$$
and with the change of variables $y=e^t$, this is 
$$
= \int_0^{\infty} \Big| \sum_{n\le e^t} \tf(n) \log n\Big|^2 e^{-2(1+\alpha)t}
dt. \tag{3.9}
$$
By Plancherel's formula
$$
\int_0^{\infty} \Big| \sum_{n\le e^t} \tf(n) \log n\Big|^2 e^{-2(1+\alpha)t}
dt = \frac{1}{2\pi} \int_{-\infty}^{\infty} 
\Big| \frac{{ F}^{\prime} (1+\alpha+iy)}{1+\alpha +iy}\Big|^2 dy. \tag{3.10}
$$

\proclaim{Lemma 3.2}  Let $T\ge 1$ be a real number.  Then 
$$
\frac{1}{2\pi} \int_{-\infty}^{\infty} \Big| \frac{F'(1+\alpha+iy)}{1+\alpha
+iy}\Big|^2 dy 
\le \Big(\max_{|y|\le T} |F(1+\alpha+iy)|^2\Big) \Big(\frac{1-x^{-2\alpha}}
{2\alpha}\Big) + O\Big(\frac 1T + \frac{m^3}{T^2} +m^2\Big)
$$
where, for convenience, we have set $m=m(\alpha)=\min(\log x,1/\alpha)$.
\endproclaim
\demo{Proof} We split the integral to be bounded into two parts: $|y|\le T$, 
and $|y|>T$.  Split the second region further into intervals of the 
form $kT \le |y|\le (k+1)T$ where $k\ge 1$ is an integer.  Thus
$$
\align
\int_{|y|>T} \Big|\frac{F'(1+\alpha+iy)}{1+\alpha+iy}\Big|^2
dy &\ll \sum_{k=1}^{\infty} \frac{1}{k^2 T^2} \int_{|y|=kT}^{(k+1)T} |F'(1
+\alpha+iy)|^2 dy \\
&\ll 
\sum_{k=1}^{\infty} \frac{1}{k^2 T^2} \sum_{n=1}^{\infty} \frac{|\tf(n)|^2
\log^2 n}{n^{2+2\alpha}} (T+n),
\\
\endalign
$$
by appealing to Corollary 3 of Montgomery and Vaughan [13].  Since $\tf(n) =0$ 
if $n$ is divisible by a prime larger than $x$, this is 
$$
\ll \frac{1}{T} \sum_{n=1}^{\infty} \frac{\log^2 n}{n^{2+2\alpha}}
+ \frac{1}{T^2} \sum\Sb n=1\\ p|n\implies p\le x\endSb^{\infty}\frac{\log^2 n}
{n^{1+2\alpha}} 
\ll \frac{1}{T} + \frac{m^3}{T^2} .
\tag{3.11}
$$

We now turn to the first region $|y|\le T$.  Define $g(n)$ to be 
the completely multiplicative function given on primes $p$ by $g(p)=\tf(p)$.
Put $G(s) =\sum_{n=1}^{\infty} g(n)n^{-s}$, and define $H(s)$ by 
$F(s)=G(s)H(s)$.  Note that $H(s)$ is absolutely convergent in 
Re$(s)>\frac 12$, and that in the region Re$(s)\ge 1$ we have uniformly
$|H(s)|$, $|H'(s)|\ll 1$.  Using $F'=G'H+GH'= F (G'/G) + O(G)$, we see that 
$$
\align
\frac{1}{2\pi} \int_{-T}^{T} \Big|\frac{F'(1+\alpha+iy)}{1+\alpha+iy}\Big|^2
dy 
&\le \Big( \max_{|y|\le T} |F(1+\alpha+iy)|^2\Big) 
\frac{1}{2\pi}\int_{-\infty}^{\infty} \Big| \frac{\frac{G'}{G}(1+\alpha+iy)}
{1+\alpha+iy}\Big|^2 dy \\
&\hskip .5 in + O\Big(\int_{-T}^{T} 
\Big|\frac{G(1+\alpha+iy)}{1+\alpha+iy}\Big|^2 dy\Big). \tag{3.12}
\\
\endalign
$$
Splitting the interval $[-T,T]$ into subintervals of length $1$, we see 
that the remainder term above is
$$
\ll \sum_{k=-[T]-1}^{[T]} \frac{1}{1+k^2} \int_{k}^{k+1} |G(1+\alpha+iy)|^2 dy
\ll \sum_{k=-[T]-1}^{[T]} \frac{1}{1+k^2} \sum_{n=1}^{\infty} 
\frac{|g(n)|^2}{n^{2+2\alpha}} (1+n)
$$
by appealing again to Corollary 3 of [13].  Plainly this is 
$$
\ll \sum\Sb n=1\\ p|n\implies p\le x\endSb^{\infty} 
\frac{1}{n^{1+2\alpha}} \ll m.  \tag{3.13}
$$

We focus on the main term in the right side of (3.12).
Since $-\frac{G'}{G}(s) \!=\! \sum_{n=1}^{\infty} g(n) \Lam(n)n^{-s}$ we 
get, by Plancherel's formula, that 
$$
\frac{1}{2\pi}\int_{-\infty}^{\infty} \Big| \frac{\frac{G'}{G}(1+\alpha+iy)}
{1+\alpha+iy}\Big|^2 dy = 
\int_0^{\infty} \Big| \sum_{n\le e^t} g(n) \Lam(n)\Big|^2 e^{-2(1+\alpha)t}
dt.
$$
Since $|g(n)|\le 1$ always, we see that $|\sum_{n\le e^t} g(n) \Lam(n)|\le 
\psi(e^t)$ for all $t$.  Further, since $g(n)=0$ if $n$ is divisible 
by a prime larger than $x$, we see that if $t\ge \log x$, then 
$|\sum_{n\le e^t} g(n)\Lam(n)| \le \psi(x)+O(\sqrt{x})$.  Using these 
observations together with the prime number theorem we deduce that 
the above is 
$$
\align
&\le \int_0^{\log x} \Big( e^t + O\Big(\frac{e^t}{(t+1)^2}\Big)\Big)^2 
e^{-2(1+\alpha) t} dt + \int_{\log x}^{\infty} \Big( x+O\Big(\frac{x}{\log^2 x}
\Big)\Big)^2 e^{-2(1+\alpha)t} dt 
\\
&= \frac{1-x^{-2\alpha}}{2\alpha}  + O(1).
\\
\endalign
$$
Thus the main term in the right side of (3.12) is 
$$
\le \Big( \max_{|y|\le T} |F(1+\alpha+iy)|^2\Big) \Big( \frac{1-x^{-2\alpha}}
{2\alpha} + O(1)\Big). \tag{3.14}
$$
Combining this with (3.13), and (3.11), we obtain the Lemma,
since $|F(1+\alpha+iy)| \ll \prod_{p\le y} (1-1/p^{1+\alpha})^{-1} \ll m$.
\enddemo

Combining (3.10) with Lemma 3.2, we conclude that (3.9) is less than 
$$
\Big( \max_{|y|\le T} |F(1+\alpha +iy)|^2 \Big)
\Big(\frac{1-x^{-2\alpha}}{2\alpha}\Big) 
+ O\Big(\frac{1}{T} 
+ \frac{m(\alpha)^3}{T^2}+m(\alpha)^2\Big). 
$$
We input this estimate into (3.8), and then use that in (3.7).  Noting that 
$(1-x^{-2\alpha})/\alpha \ll m(\alpha)$,
we conclude that 
$$
\align
\int_2^x \frac{|S(y)|}{y^2} dy &\le 2 \int_0^1 
\Big( \max_{|y|\le T} |F(1+\alpha +iy)| \Big)
\Big(\frac{1-x^{-2\alpha}}{2\alpha}\Big) d\alpha 
\\
&\qquad \qquad + O\Big(\int_0^{1} \Big(\frac{\sqrt{m(\alpha)}}{\sqrt T} 
+ \frac{m(\alpha)^2}{T} + m(\alpha) +\log\log x \Big) d\alpha\Big) \\
&=2 \int_0^1 
\Big( \max_{|y|\le T} |F(1+\alpha +iy)| \Big)
\Big(\frac{1-x^{-2\alpha}}{2\alpha}\Big) d\alpha +
O\left(  \frac{\log x}T + \log\log x \right) .
\\
\endalign
$$
When used with (2.1) of Lemma 2.1, this yields Proposition 1.

We end this section by giving a variant of Propostion 1 which will be 
our main tool in the proof of Theorem 4.

\proclaim{Proposition 3.3}  Let $f$, $T$, and $x$ be as in Theorem 1.
Then for $1\le w \le x$, we have
$$
\align
\Big| \frac{S(x)}{x} -\frac{S(x/w)}{x/w}\Big| 
&\ll \frac{1}{\log x} \int_{0}^1 m(\alpha)  
\Big( \max_{|y|\le T} |(1-w^{-\alpha-iy}) F(1+\alpha+iy)| \Big)d\alpha \\
&\hskip .5 in
+  O\Big( \frac{1}{T}+ \frac{\log 2w}{\log x} \log \frac{\log x}{\log 2w}\Big).
\\
\endalign
$$
\endproclaim
\demo{Proof}  Since the proof is very similar to that of Proposition 1, we 
shall merely sketch it.  Arguing as in (3.7), we get that 
$$
\align
\int_{2w}^{x} \Big|\frac{S(y)}{y}-\frac{S(y/w)}{y/w} \Big|
\frac{dy}{y}
&\ll \int_0^1 \Big( \int_{2w}^x \Big| \frac{1}{y} \sum_{n\le y} f(n)\log n
-\frac{1}{y/w} \sum_{n\le y/w} f(n)\log n \Big| \frac{dy}{y^{1+2\alpha}}
\Big) d\alpha \\
&\hskip .5 in +  \log 2w \log \Big(\frac{\log x}{\log 2w}\Big).
\\
\endalign
$$
Using Cauchy's inequality as in (3.8), we see that 
$$
\align
&\int_{2w}^x \Big| \frac{1}{y} \sum_{n\le y} f(n)\log n
-\frac{1}{y/w} \sum_{n\le y/w} f(n)\log n \Big| \frac{dy}{y^{1+2\alpha}}
\\
\ll &\min m(\alpha)^{\frac 12} 
\Big( \int_{2w}^{x} \Big| \frac{1}{y} \sum_{n\le y} f(n)\log n
-\frac{1}{y/w} \sum_{n\le y/w} f(n)\log n \Big|^2 \frac{dy}{y^{1+2\alpha}}
\Big)^{\frac {1}{2}}.\\
\endalign
$$
As before, we handle the second factor above by replacing $f$ by $\tf$,
extending the range of integration to $\int_1^{\infty}$, substituting $y=e^t$, 
and invoking Plancherel's formula.  The only difference from (3.10) is that 
$F'(1+\alpha+iy)/(1+\alpha+iy)$ in the right side there 
must be replaced by the Fourier transform of 
$e^{-(1+\alpha)t} \sum_{n\le e^t} \tf(n) \log n - we^{-(1+\alpha)t} \sum_{n\le 
e^t/w} \tf(n) \log n$ which is 
$- F'(1+\alpha+iy) (1-w^{-\alpha-iy})/(1+\alpha+iy)$.  We make this 
adjustment, and follow the remainder of the proof of Proposition 1.

\enddemo

\head 4. Proofs of Theorem 1 and Corollary 1\endhead

\noindent Recall the multiplicative function $\tf(n)$ defined 
by $\tf(p^k)=f(p^k)$ for $p\le x$, and  $\tf(p^k)=0 $ for $p>x$.  Then 
$F(s) =\sum_{n} \tf(n)n^{-s}$, and since $|\tf(n)|\le 1$ always, 
we get that for all $0 < \alpha \le 1$, 
$$
\max_{y \in {\Bbb R}} |F(1+\alpha+iy)| 
\le \zeta(1+\alpha) = \frac{1}{\alpha}+ O(1). \tag{4.1}
$$
Taking $a_n=\tf(n)$ in Lemma 2.2 and noting that $\sum_{n} |a_n|/n 
\ll \log x$, we conclude that for $0\le \alpha \le 1$
$$
\max_{|y| \le T} |F(1+\alpha+iy)| \le 
\max_{|y|\le 2T} |F(1+iy)| + O\Big(\frac{\alpha \log x}{T}\Big).
\tag{4.2}
$$

Note that $L\le \frac{1}{\log x} \prod_{p\le x}(1-\frac 1p)^{-1} = e^{\gamma} 
+O(1/\log x)$, by Mertens' theorem.  The Theorem is trivial if 
$1\le L \le e^{\gamma}+O(1/\log x)$, and we suppose that $L\le 1$.  
We use Proposition 1, 
employing the bound (4.2) when $\alpha \le 1/(L\log x)$, and the 
bound (4.1) when $1/(L\log x) \le \alpha \le 1$.  We deduce that 
$$
\frac{1}{x} \Big|\sum_{n\le x} f(n) \Big|  \le L
\int_0^{1/L\log x}  
\frac{1-x^{-2\alpha}}{\alpha} d\alpha 
+ \frac{2}{\log x} \int_{1/L\log x}^1 \frac{1-x^{-2\alpha}}{2\alpha} \frac {1}
{\alpha} d\alpha  + O\Big(\frac{1}{T} 
+ \frac{\log \log x}{\log x}\Big).
\tag{4.3}
$$

Making a change of variables $y=2\alpha \log x$, we see that the 
first integral above is 
$$
\align
\le L \int_0^{2/L} \frac{1-e^{-y}}{y} dy 
&= L
\Big( \int_0^1 \frac{1-e^{-y}}{y} dy + \int_1^{2/L} \frac{dy}{y} 
 - \int_1^{\infty} \frac{e^{-y}}{y} dy + \int_{2/L}^{\infty}
 \frac{e^{-y}}{y}
dy\Big)\\
&= L \Big(\gamma+\log \frac{2}{L}\Big) + L \int_{2/L}^{\infty} 
\frac{e^{-y}}{y} dy,\\
\endalign 
$$
since $\gamma = \int_0^1 (1-e^{-y})/y dy - \int_1^{\infty} e^{-y}/y dy$.  
Further, the second integral in (4.3) is 
$$
\align
&\frac{1}{\log x}
\int_{1/L\log x}^1 \frac{1-x^{-2\alpha}}{\alpha^2} d\alpha
= 2\int_{2/L}^{\infty} \frac{1-e^{-y}}{y^2} dy 
=  L - \int_{2/L}^{\infty} \frac{2e^{-y}}{y^2} dy.\\
\endalign
$$
Combining the above bounds, we see that the right side of (4.3) is  
$$
\leq L\Big(1+\log 2 + \log \frac{e^{\gamma}}{L} 
+ \int_{2/L}^{\infty} \frac{e^{-y}}{y} 
\Big( 1 - \frac{2/L}{y} \Big) dy \Big)
+ O\Big(\frac{1}{T}+\frac{\log\log x}
{\log x}\Big).
\tag{4.4}
$$
Since the maximum of $(1-(2/L)/y)/y$ for $y\ge 2/L$ is 
attained at $y=4/L$, we see that the integral term above is 
$\le L/8 \int_{2/L}^{\infty} e^{-y} dy \leq 
1/8 \int_{2}^{\infty} e^{-y} dy = 1/(8e^2)$ since $L\le 1$, and theorem
then follows from (4.4) since $1+\log 2 + 1/(8e^2) \le 12/7$.

We now deduce Corollary 1.  Suppose $f$ is completely multiplicative.  
Then, by Mertens' theorem, 
$$
\align
|F(1+iy)| &= (e^{\gamma}\log x +O(1)) \prod_{p\le x} 
\Big|1-\frac{f(p)}{p^{1+iy}}\Big|^{-1}\Big(1-\frac 1p\Big) 
\\
&= (e^{\gamma} \log x +O(1)) \exp\Big(- \sum\Sb p \le x\\ k\ge 1\endSb 
\frac{1-\RE f(p^k)p^{-iky}}{kp^k} \Big),\tag{4.5}\\ 
\endalign
$$
and so it follows that $L\le e^{\gamma -M} + O(1/\log x)$.  
Using this bound in Theorem 1, we get the completely multiplicative case of 
Corollary 1.

If $f$ is only known to be multiplicative then note that 
$$
\Big| 1+\frac{f(p)}{p^{1+iy}} + \frac{f(p^2)}{p^{2+2iy}} +\ldots \Big|
\Big| 1-\frac{f(p)}{p^{1+iy}}\Big| 
\le 1 + \frac{2}{p(p-1)},
$$
since $|f(p^k)|\le 1$ for all $k$.  Using this with the observation 
of the preceding paragraph, we see that $L\le \prod_p (1+\frac{2}{p(p-1)}) 
e^{\gamma -M} +O(1/\log x)$ in this case.  
Appealing now to Theorem 1, and noting that $\log (\prod_{p} 
(1+\frac{2}{p(p-1)}))
\ge 8/7$, we deduce this case of Corollary 1.

\head 5. Proof of Theorem 3 \endhead

\noindent We may suppose that $|y_0| \ge 10$.  Applying Theorem 1 with 
$T=|y_0|/2 -1$ we get that 
$$
\align
\frac{1}{x} \Big| \sum_{n\le x} f(n) \Big| 
&\ll \Big( \frac{\max_{|y|\le |y_0|-2} |F(1+iy)|}{\log x}\Big) 
\log \Big(\frac{e^{1+\gamma}\log x}{\max_{|y|\le y_0-2} |F(1+iy)|} \Big) 
\\
&\hskip .5 in + \frac{1}{|y_0| +1} + \frac{\log \log x}{\log x}. 
\tag{5.1}\\
\endalign
$$
By the definition of $y_0$, we see that for $|y|\le |y_0|-2$, 
$$
|F(1+iy)| \le \Big( |F(1+iy)F(1+iy_0)|\Big)^{\frac{1}{2}},
$$
and appealing to Lemma 2.3, this is (with $\log x \gg |\beta|=|y-y_0| \ge 2$)
$$
\ll (\log x)^{\frac{2}\pi} (\log \log x)^{2(1-\frac{2}{\pi})}.
$$ 
Using this bound in (5.1), we obtain the Theorem.

\head 6. Proof of Theorem 4 \endhead

\noindent  If $|y_0| \ge (\log x)/2$, then 
in view of Theorem 3, the result follows.  Thus we may assume that 
$|y_0|\le (\log x)/2$.  Put $f_{0}(n) = f(n)n^{-iy_0}$, and 
define $F_{0}(s) = \prod_{p\le x} (1+f_{0}(p)p^{-s} + f_{0}(p^2)p^{-2s} 
+\ldots) = F(s+iy_0)$.  We note that 
$$
(|F(1+iy_0)| =) |F_0(1)| = \max_{|y|\le \log x} |F_0(1+iy)|. \tag{6.1}
$$
Indeed, the left side of (6.1) is plainly $\le $ right side; and further 
the right side is $= \max_{|y|\le \log x} |F(1+iy+iy_0)| \le 
\max_{|y|\le \log x +|y_0|} |F(1+iy)| \le |F(1+iy_0)|$, proving (6.1).

We now appeal to Proposition 3.3, with $f$ there replaced by $f_0$, and 
$F$ by $F_0$, and with $T= (\log x)/2$.  Thus we see that 
$$
\align
\Big|\frac{1}{x} \sum_{n\le x} f_0(n) - &\frac{w}{x} 
\sum_{n\le x/w} f_0(n)\Big | 
\ll  \frac{\log 2w}{\log x} \log \Big(\frac{\log x}{\log 2w}\Big)\\
&+\frac{1}{\log x} \int_0^1 \min\Big(\log x,\frac{1}{\alpha}\Big) 
\Big(\max_{|y|\le (\log x)/2} |F_0(1+\alpha+iy) (1-w^{-\alpha-iy})|\Big) 
d\alpha. 
\tag{6.2}\\
\endalign
$$

Next, we use Lemma 2.2 with $a_n=f_0(n)$ if $n$ is divisible 
only by primes $\le x$, and $a_n=0$ otherwise.  Thus $A(s)=F_0(s)$, 
and $\sum_{n=1}^{\infty} |a_n|/n \ll \log x$.  Taking $T=(\log x)/2$, 
we deduce from (2.6) of Lemma 2.2 that 
$$
\max_{|y|\le (\log x)/2} |F_0(1+\alpha+iy) (1-w^{-\alpha-iy})| 
\le \max_{|y|\le \log x} |F_0(1+iy) (1-w^{-iy})|+ O(1). \tag{6.3}
$$
If $|y|\le 1/\log x$, then plainly $|F_0(1+iy)(1-w^{-iy})| \ll 
\log x (|y|\log 2w) \ll \log 2w$.  If $\log x \ge |y|> 1/\log x$, then 
using (6.1) and Lemma 2.3, we get 
$$
|F_0(1+iy)| \le \Big( |F_0(1)F_0(1+iy)|\Big)^{\frac{1}{2}} 
\ll (\log x)^{\frac 2\pi} \max\Big(\frac{1}{|y|}, (\log \log x)^{2}
\Big)^{(1-\frac{2}{\pi})}.
$$
Since $|1-w^{-iy}| \ll \min(1,|y|\log 2w)$, we deduce from these 
remarks and (6.3) that 
$$
\max_{|y|\le (\log x)/2} |F_0(1+\alpha+iy) (1-w^{-\alpha-iy})| 
\ll (\log x)^{\frac{2}{\pi}} \max\Big(\log 2w, (\log \log x)^2
\Big)^{1-\frac{2}{\pi}}.
\tag{6.4}
$$

In addition, we have the trivial estimate 
$$
\max_{|y|\le (\log x)/2} |F_0(1+\alpha+iy) (1-w^{-\alpha-iy})| \ll 
\zeta(1+\alpha) \ll \frac{1}{\alpha}. \tag{6.5}
$$ 
We now use (6.2), employing estimate (6.4) when 
$\alpha$ is less than 

\noindent $\max\Big(\log 2w,(\log \log x)^2\Big)^{-(1-\frac{2}{\pi})} 
(\log x)^{-\frac{2}{\pi}}
$, and estimate 
(6.5) for larger $\alpha$.  This gives the Theorem.

\head 7. Deduction of Corollary 3 \endhead

\noindent We require the following lemma, which relates the mean value of 
$f(n)$
to the mean-value of $f(n)n^{i\alpha}$.

\proclaim{Lemma 7.1}  Suppose $f(n)$ is  a multiplicative function with
$|f(n)| \le 1$ for all $n$.  Then for any real number $\alpha$ we have
$$
\sum_{n\le x} f(n) n^{i\alpha} = 
\frac{x^{i\alpha}}{1+i\alpha} \sum_{n\le x } f(n) + O\biggl(\frac{x}{\log x}
 \log (e+|\alpha|)
\exp\biggl(\sum_{p\le x} \frac{|1-f(p)|}{p}\biggr)\biggr).
$$
\endproclaim

To prove this Lemma, we require a consequence of Theorem 2 of 
Halberstam and Richert [7].  Suppose $h$ is a non-negative 
multiplicative function with $h(p^k)\leq 2$ for all prime powers $p^k$. 
It follows from Theorem 2 of [7] that 
$$
\sum_{n\leq x} h(n) \leq 
\ \frac{2x}{\log x} \ \sum_{n\leq x} \frac{h(n)}{n} \ 
\left\{ 1 + O\left( \frac{1}{\log x}\right) \right\}. \tag{7.1}
$$
Using partial summation we deduce from (7.1) that for $1\leq y\leq x^{1/2}$,
$$
\sum_{x/y< n\leq x} \frac{h(n)}{n} \leq \left\{ \frac{1}{\log x} - 
\log \left( 1 - \frac{\log y}{\log x} \right) \right\}
 \ \sum_{n\leq x} \frac{h(n)}{n}  
\left\{ 2 + O\left( \frac{1}{\log x}\right) \right\} .\tag{7.2}
$$

\demo{Proof of Lemma 7.1}  
Let $g$ denote the multiplicative function defined by $g(p^k) = 
f(p^k) -f(p^{k-1})$, so that $f(n) =\sum_{d|n} g(d)$.  Then
$$
\sum_{n\le x} f(n) n^{i\alpha} = \sum_{n\le x} n^{i\alpha} \sum_{d|n} g(d)
=\sum_{d\le x} g(d)d^{i\alpha} \sum_{n\le x/d} n^{i\alpha}.  \tag{7.3}
$$
By partial summation it is easy to see that
$$
\sum_{n\le z }n^{i\alpha} = 
\cases
\frac{z^{1+i\alpha}}{1+i\alpha} +O(1+\alpha^2)\\
O(z).\\
\endcases
$$
We use the first estimate above in (7.3) when $d\le x/(1+\alpha^2)$, 
and the second estimate when $x/(1+\alpha^2)\le d\le x$.  
This gives
$$
\align
\sum_{n\le x} f(n)n^{i\alpha} 
&= \frac{x^{1+i\alpha}}{1+i\alpha} \sum_{d\le x}\frac{g(d)}{d} 
+ O\biggl((1+\alpha^2)\sum_{d\le x/(1+\alpha^2)} |g(d)| + 
x\sum_{x/(1+\alpha^2)\le d\le x} 
\frac{|g(d)|}{d}\biggr).\\
\endalign
$$

Applying (7.1) and (7.2) we deduce that
$$
\align
\sum_{n\le x} f(n)n^{i\alpha} &=\frac{x^{1+i\alpha}}{1+i\alpha} \sum_{d\le x} 
\frac{g(d)}{d} + O\biggl(\frac{x}{\log x}\log (e+|\alpha|) \sum_{d\le x} 
\frac{|g(d)|}{d}\biggr)\\
&=\frac{x^{1+i\alpha}}{1+i\alpha} \sum_{d\le x} \frac{g(d)}{d} + O\biggl(
\frac{x}{\log x} \log (e+|\alpha|) \exp\biggl(\sum_{p\le x} \frac{|1-f(p)|}{p}
\biggr)\biggr).\\
\endalign
$$

Using the above estimate twice, once with $\alpha$ replaced by $0$, we 
obtain the Lemma.

\enddemo

We now proceed to the proof of Corollary 3.  We may suppose that 
$w\le \sqrt{x}$, else there's nothing to prove.  Let $y_0$ be as 
in Theorem 4. By the definition of $M$ and by (4.5)
 we know that for all $|y|\le 2\log x$, 
$$
\sum_{p\le x} \frac{1-\RE f(p) p^{-iy}}{p} \ge M = \sum_{p\le x} 
\frac{1-\RE f(p)p^{-iy_0}}{p} + O(1).
$$
Further we have for $|y|\le 2\log x$
$$
\sum_{p\le x/w} \frac{1-\RE f(p)p^{-iy}}{p} \ge 
\sum_{p\le x} \frac{1-\RE f(p)p^{-iy}}{p} - 2\sum_{x/w\le p\le x} \frac{1}p
\ge M + O(1).
$$
By Corollary 1 (with $T=\log x$) it follows that 
$$
\frac{1}{x} \Big| \sum_{n\le x} f(n) \Big|, \qquad \frac{1}{x/w} 
\Big|\sum_{n\le x/w} f(n)\Big|  \ll Me^{-M} + 
\frac{\log \log x}{\log x}.
$$
 From this estimate, Corollary 3 follows if $M\ge (2-\sqrt{3}) \log \log x$.  
We suppose now that $M\le (2-\sqrt{3}) \log \log x$.  

For a complex number $z$ in the unit disc, we have $|1-z| = 
(1+|z|^2-2\RE z)^{\frac{1}{2}} \le (2-2\RE z)^{\frac12}$.  Hence, 
by Cauchy's inequality and our bound on $M$, 
$$
\align
\sum_{p\le x} \frac{|1-f(p)p^{-iy_0}|}{p}  &\le \sum_{p\le x} 
\frac{\sqrt{2-2\RE f(p)p^{-iy_0}}}{p} 
\le \Big( \sum_{p\le x} \frac{2}{p} \Big)^{\frac 12} 
\Big( \sum_{p\le x} \frac{1-\RE f(p)p^{-iy_0}}{p} \Big)^{\frac 12}\\ 
&\le \Big(2(2-\sqrt{3})\Big)^{\frac 12} \log \log x + O(1)= 
(\sqrt{3}-1) \log \log x + O(1).
\tag{7.4}\\
\endalign
$$ 
Applying Lemma 7.1, we see that 
$$
\align
\frac{1}{x} \sum_{n\le x} f(n) &= \frac{x^{iy_0}}{1+iy_0} \sum_{n\le x} f(n) 
n^{-iy_0} + O\Big( \frac{\log \log x}{\log x} \exp\Big( \sum_{p\le x} 
\frac{|1-f(p)p^{-iy_0}|}{p} \Big)\Big)\\ 
&= \frac{x^{iy_0}}{1+iy_0} \sum_{n\le x} f(n) 
n^{-iy_0} + O\Big( \frac{\log \log x}{(\log x)^{2-\sqrt{3}}}\Big),
\\
\endalign
$$
and similarly
$$
\frac wx \sum_{n\le x/w} f(n) = \frac{(x/w)^{iy_0}}{1+iy_0} 
\sum_{n\le x/w} f(n)n^{-iy_0} + O\Big(\frac{\log \log x}
{(\log x)^{2-\sqrt{3}}}\Big).
$$
Taking absolute values in these relations, and appealing to Theorem 4, we 
obtain the Corollary.

\head 8. Proof of Theorem 5 \endhead 

\noindent We recall the notations of \S 3a.  
We first obtain a lower bound for $M(t)$ in terms 
of $M_0= \int_0^u \frac{1-\RE \chi(v)}{v} dv$.

\proclaim{Proposition 8.1}  For all $t>0$ we have 
$$
M(t) \ge \max \Big( 0, \kappa M_0 -\kappa \nu \log(tu) +(1-\kappa \nu)
\int_{tu}^{\infty} \frac{e^{-v}}{v} dv + C(D)\Big),
$$
where $C(D)$ was defined in (1.6).
\endproclaim

\demo{Proof} First note that $M(t)\geq 0$ by definition. Also 
$$
M(t) - \kappa M_0 = I - \kappa \int_0^u 
\frac{1-\RE \chi(v)}{v} (1-e^{-tv}) dv + \int_u^{\infty} \frac{e^{-tv}}{v} 
dv
$$
where 
$$
I:= \min_{y\in {\Bbb R}} \int_0^u \frac{1-\kappa - \RE \chi(v) 
(e^{-ivy}-\kappa)}{v} e^{-tv} dv. \tag{8.1}
$$
Since $1-\RE \chi(v) \le \nu$, we get that 
$$
\align 
-\kappa \int_0^u &\frac{1-\RE \chi(v)}{v} (1-e^{-tv}) dv + \int_u^{\infty} 
\frac{e^{-tv}}{v} dv \\
&\ge -\kappa \nu \Big(\int_0^u \frac{1-e^{-tv}}{v} dv -\int_u^{\infty} 
\frac{e^{-tv}}{v} dv \Big) + (1-\kappa \nu) \int_u^{\infty} \frac{e^{-tv}}{v} 
dv\\
&= -\kappa \nu (\gamma +\log (tu)) + (1-\kappa \nu) \int_{tu}^{\infty} 
\frac{e^{-v}}{v} dv,
\\
\endalign
$$
so that 
$$
M(t) \ge \kappa M_0 - \kappa \nu (\gamma +\log (tu)) + (1-\kappa \nu) 
\int_{tu}^{\infty} \frac{e^{-v}}{v} dv + I.
$$
Therefore we obtain the Proposition by proving 
$$
I \ge \min_{\epsilon = \pm 1}
\int_0^{2\pi} \frac{\min(0,1-\kappa -\max_{\delta \in D} 
\RE \delta(e^{\epsilon ix} -\kappa))}{x} dx. \tag{8.2}
$$

If the minimum in (8.1) occurs at $y=0$, then 
$I =(1-\kappa) \int_0^u \frac{ 1-\RE \chi(v)}{v} e^{-tv} dv \ge 0$,
which is stronger than (8.2).  So we may suppose that the 
minimum in (8.1) occurs for some $y\neq 0$.  Put $w(\theta) = 1-\kappa 
- \max_{\delta\in D} \RE \delta(e^{-i\theta}-\kappa)$.  Then we see that,
with $\epsilon = \text{sgn}(y)$,
$$
I \ge \int_0^u \frac{w(vy)}{v} e^{-tv} dv = \int_0^{|y|u} 
\frac{w(\epsilon v)}{v} e^{-tv/|y|} dv = \int_0^{|y|u} \frac {e^{-tv/|y|}}{v} 
d\Big(\int_0^v w(\epsilon x) dx\Big).
$$
Integrating by parts, we conclude that 
$$
I \ge \frac{e^{-tu}}{|y|u} \int_0^{|y|u} w(\epsilon x)dx + \int_0^{|y|u} 
\Big(\int_0^v w(\epsilon x)dx \Big) \Big( \frac{e^{-tv/|y|}}{v^2} + 
\frac{te^{-tv/|y|}}{v|y|}\Big) dv. \tag{8.3}
$$

Note that $w(\epsilon x)$ is a $2\pi$-periodic function, and that 
$\frac{1}{2\pi} \int_0^{2\pi} w(\epsilon x) dx = 1- \hbar(\kappa) \ge 0$.  
Hence putting $w^{-} (x) =\min (0,w(x))$, we get that 
$$
\int_0^v w(\epsilon x) dx \ge \int_{2\pi [\frac{v}{2\pi}]}^v w^{-}(\epsilon x)
dx = W_{\epsilon} (v), \tag{8.4}
$$
say.  Observe that $W_{\epsilon}$ is a $2\pi$-periodic function, which is 
always negative, and that $W_{\epsilon}$ is decreasing in $(0,2\pi)$.  

Using (8.4) in (8.3), and since $W_\epsilon$ is negative and $e^{-x}(1+x) 
\le 1$ for all $x\ge 0$ , we get that
$$
I \ge \frac{e^{-tu}}{|y|u} W_\epsilon(|y|u) +\int_0^{|y|u} 
\frac{W_{\epsilon}(v)}{v^2} e^{-tv/|y|} \Big(1+\frac{tv}{|y|}\Big) dv 
\ge \frac{W_\epsilon(|y|u)}{|y|u} +\int_0^{|y|u} \frac{W_\epsilon (v)}{v^2}dv.
\tag{8.5}
$$

If $\alpha \ge 2\pi$ then, since $W_\epsilon (v) \ge W_\epsilon(2\pi-)$
($=\int_0^{2\pi} w^{-}(\epsilon x)  dx$),
we get 
$$
\align
\frac{W_\epsilon(\alpha)}{\alpha} +\int_0^\alpha \frac{W_\epsilon(v)}{v^2}dv
&\ge 
\frac{W_\epsilon(2\pi-)}{\alpha} + \int_0^{2\pi} \frac{W_\epsilon(v)}{v^2} 
dv + W_\epsilon(2\pi -) \int_{2\pi}^{\alpha} \frac{dv}{v^2} 
\\
&= \int_0^{2\pi} \frac{W_\epsilon(v)}{v^2} + \frac{W_\epsilon(2\pi-)}{2\pi}.\\
\endalign
$$
If $\alpha <2\pi$ then, since $W_\epsilon(x)$ is decreasing in $(0,2\pi)$,
$$
\align
\!
\frac{W_\epsilon (\alpha)}{\alpha} +\int_0^\alpha \frac{W_\epsilon(v)}{v^2} dv 
&\ge \frac{W_\epsilon(\alpha)}{\alpha} +\int_0^{2\pi} \frac{W_\epsilon(v)}{v^2}
- W_\epsilon(\alpha) \int_\alpha^{2\pi} \frac{dv}{v^2 } 
= \int_0^{2\pi} \frac{W_\epsilon(v)}{v^2} dv +\frac{W_\epsilon(\alpha)}{2\pi}
\!\\
&\ge \int_0^{2\pi} \frac{W_\epsilon(v)}{v^2} dv + \frac{W_\epsilon(2\pi-)}
{2\pi}.\\
\endalign
$$
Using these in (8.5), we conclude that 
$$
I \ge \int_0^{2\pi} \frac{W_\epsilon(v)}{v^2} dv+ 
\frac{W_\epsilon(2\pi-)}{2\pi}
= \int_0^{2\pi} \frac{w^-(\epsilon x)}{x} dx,
$$
which, from the definition of $w^-$, is greater than or equal to the right 
side of (8.2) for both $\epsilon =\pm 1$.  This 
completes the proof of the Proposition.  

\enddemo

We now finish the proof of Theorem 5.  We first deal with the 
case $D\neq [0,1]$, where $\kappa \nu <1$.  
We shall input the bounds for $M(t)$
in  Proposition 8.1 into the $t$-integral in Proposition 3.1.  We split this 
integral into three parts: when $0 \le t\le t_1:=e^{-\gamma}/u$, 
when $t_1 \le t\le t_2:=\exp(\frac{M_0}{\nu}+\frac{C(D)}{\kappa\nu})
/u$, and when $t> t_2$.

We first estimate the contribution of the first range of $t$.  Since 
$$
\int_{tu}^{\infty} \frac{e^{-v}}{v} dv \ge \int_{tu}^{1} \frac{dv}{v} 
- \int_{0}^{1} \frac{1-e^{-v}}{v} dv +\int_1^{\infty} \frac{e^{-v}}{v} dv
= -\log (tu) -\gamma,
$$
and $1-\kappa \nu \ge 0$, we see that 
$M(t) \ge \kappa M_0 - \log (tu) +C(D) -\gamma(1-\kappa \nu)$, by Proposition
8.1.  Hence, with a little calculation,
$$
\align
\int_0^{t_1} \Big( \frac{1-e^{-2tu}}{t} \Big) \frac{\exp(-M(t))}{tu}
dt &\le \exp(-\kappa M_0 -C(D) +\gamma(1-\kappa \nu)) \int_0^{2e^{-\gamma}} 
\frac{1-e^{-x}}{x} dx \\
&\le \exp(-\kappa M_0 -C(D) +\gamma(1-\kappa \nu)).
\\
\endalign
$$

For the middle range of $t$, we use the bound $M(t)\ge \kappa M_0 -\kappa \nu 
\log(tu) +C(D) $, which holds since $1-\kappa \nu$, and $\int_{tu}^{\infty}
\frac{e^{-v}}{v} dv$ are non-negative.  Hence
$$
\align
\int_{t_1}^{t_2}
\Big(\frac{1-e^{-2tu}}{t}\Big) \frac{\exp(-M(t))}{tu} dt 
&\le \exp(-\kappa M_0 -C(D) )
\int_{t_1}^{t_2}  \frac{(tu)^{\kappa \nu}}{tu}\frac{dt}{t}\\
&=\frac{\exp(-\kappa M_0-C(D)+\gamma(1-\kappa\nu))}{1-\kappa\nu}
- \frac{\exp(-\frac{M_0}{\nu}-\frac{C(D)}{\kappa \nu})}{1-\kappa\nu}.
\\
\endalign
$$

For the last range of $t$, we use the trivial bound $M(t)\ge 0$.  This 
gives that 
$$
\int_{t_2}^\infty
\Big(\frac{1-e^{-2tu}}{t}\Big) \frac{\exp(-M(t))}{tu} dt
\le \int_{t_2}^\infty
\frac{dt}{t^2 u} = \exp\Big(-\frac{M_0}{\nu}-\frac{C(D)}{\kappa\nu}\Big)
.
$$
Combining the above three bounds with Proposition 3.1, we obtain Theorem 
5 in the case $\kappa \nu <1$.

We now consider the case $D=[0,1]$ where we shall show that $|\sigma(u)|
\le e^{\gamma -M_0}$.  Put $\hchi(t)=\chi(t)$ if $t\le u$, and $\hchi(t)=0$ 
for $t>u$, and let $\hsigma$ denote the corresponding solution to (1.8).  
Note that both $\sigma(v)$ and $\hsigma(v)$ are non-negative for all $v$, 
and that $\hsigma(v)=\sigma(v)$ for $v\le u$.  
Now, using (3.6), 
$$
\align
\sigma(u) &= \frac{1}{u} \int_0^u \sigma(v) \chi(u-v) dv \le \frac 1u \int_0^u 
\sigma(v) dv \le \frac{1}{u} \int_0^{\infty} \hsigma(v)dv = 
\frac{1}{u} \lim_{t \to 0} \lap(\hsigma,t) \\
&=\frac 1u \lim_{t\to 0} 
\frac 1t \exp\Big( -\lap\Big(\frac{1-\hchi(v)}{v}, t\Big)\Big) 
= e^{-M_0} \lim_{t\to 0} \frac{1}{tu} \exp\Big(-\int_u^{\infty} 
\frac{e^{-tv}}{v} dv \Big)\\
&= e^{-M_0} \lim_{y\to 0} \frac{1}{y} \exp\Big(-
\int_y^\infty \frac{e^{-v}}{v} dv\Big) = e^{\gamma -M_0},
\\
\endalign
$$
which proves the Theorem in this case.

\head 9.  Deduction of Theorem 2 \endhead

\noindent Let $y=\exp((\log x)^{\frac 23})$, and let $g$ be the completely 
multiplicative function with $g(p)=1$ for $p\le y$, and $g(p)=f(p)$ 
for larger $p$.  Let $\chi(t)=1$ for $t\le 1$, and put for $t >1$
$$
\chi(t) = \frac{1}{\vartheta(y^t)} \sum_{p\le y^t} g(p)\log p. 
$$
Let $\sigma$ denote the corresponding solution to (1.8).  Note that 
for $u\ge 1$ 
$$
\align
\int_0^{u} \frac{1-\RE \chi(v)}{v} dv 
&= \int_1^u \frac{1}{v\vartheta(y^v)} \sum_{y < p\le y^v} (1-\RE f(p)) \log p 
\ dv  \\
&= \int_1^u \frac{1}{vy^v} \sum_{y<p \le y^v} (1-\RE f(p))\log p \ dv + 
O\Big(\frac{1}{\log y}\Big),\\
\endalign
$$ 
upon using the prime number theorem.  Interchanging the sum and the 
integral, the above is 
$$
\align
&= \sum_{y < p \le y^u} (1-\RE f(p)) \log p \int_{\log p/\log y}^u 
\frac{dv}{vy^v} + O\Big( \frac{1}{\log y}\Big) 
\\
&= \sum_{y < p \le y^u} (1-\RE f(p)) \log p \Big( 
\frac{1}{p\log p} + O\Big(\frac{1}{p\log^2 p} +\frac{1}{uy^u \log y}\Big) 
\Big) + O\Big(\frac{1}{\log y}\Big). 
\\
\endalign
$$
We conclude that 
$$
\int_0^{u} \frac{1-\RE \chi(v)}{v} dv = \sum_{y<p \le y^u} \frac{1-\RE f(p)}{p}
+ O\Big(\frac{1}{\log y}\Big). \tag{9.1}
$$

Appealing to Propositions 3 and then 2 we obtain that 
$$
\align
\frac{1}{x} \sum_{n\le x} f(n) &= \Theta(f,y) \frac{1}{x} \sum_{n\le x}
g(n)  + O\Big( \frac{1}
{(\log x)^{\frac{1}{3}}}\exp\Big(\sum_{p\le x} \frac{|1-f(p)|}{p}\Big)\Big)\\
&=\Theta(f,y) \sigma\Big(\frac{\log x}{\log y}\Big) + 
O\Big( \frac{1}
{(\log x)^{\frac{1}{3}}}\exp\Big(\sum_{p\le x} \frac{|1-f(p)|}{p}\Big)\Big).\\
\endalign
$$
Since $f(p) \in D$ for all $p$ and $D$ is convex, thus 
$\chi(t) \in D$ for all $t$.  Hence using Theorem 5 and (9.1), 
we conclude that
$$
\align
\frac{1}{x} \Big | \sum_{n\le x} f(n) \Big| 
&\le |\Theta(f,y)| \Big( \frac{2-\kappa \nu}{1-\kappa\nu} \Big) 
\exp\Big(-\kappa \sum_{y\le p\le x} \frac{1-\RE f(p)}{p} 
-C(D) +\gamma (1-\kappa \nu) \Big) \\
&\hskip .75 in +O\Big( \frac{1}
{(\log x)^{\frac{1}{3}}}\exp\Big(\sum_{p\le x} \frac{|1-f(p)|}{p}\Big)\Big). 
\\
\endalign
$$
Arguing as in (7.4), we see that 
$$
\sum_{p\le x} \frac{|1-f(p)|}p 
\le \Big( 2\log \log x \sum_{p\le x} 
\frac{1-\RE f(p)}{p}\Big)^{\frac 12} + O(1),
$$
and this completes the proof of the first part of Theorem 2.

In fact the second part of Theorem 2 follows from Lemma 2.1 for, from (2.1)
we have
$$
(\log x +1) \frac 1x \sum_{n\leq x} f(n) \leq \sum_{n\leq x} \frac{f(n)}{n}
+O(1) = e^\gamma \log x \ \Theta (f,x) + O(1)
$$
using Mertens' theorem, and the result follows.

\head 10.  Explicit Constructions \endhead

\subhead 10a.
Examining proper subregions $D$ of ${\Bbb U}$  is necessary\endsubhead

As we remarked after Theorem 5, it is not, {\sl a priori}, clear that 
one should look at proper subsets $D$ of ${\Bbb U}$ when looking
for bounds (of the shape of Theorem 2) on solutions to (1.8). However
if we   take $\chi(t) = 1$ for $t\le 1$, and $\chi(t) = e^{i\alpha t}$ for 
$t>1$ then
$$
\lap \left( \frac{1-\RE \  \chi(v) e^{-iv\alpha}}{v},t \right)
= \int_0^1 \left( \frac{1-\cos v\alpha}{v}\right) e^{-tv} dv \to
\int_0^1\frac{1-\cos v\alpha}{v}  dv \gg_\alpha 1,
$$
as $t\to 0$. Therefore, by (3.6), 
$$
\int_0^\infty |\sigma(v)| e^{-tv} dv \geq
|\lap ( \sigma, t+i\alpha)| \gg_\alpha 1/t,
$$
if $t$ is sufficiently small.  Now
$\int_{b/t}^\infty |\sigma(v)| e^{-tv} dv \leq \int_{b/t}^\infty e^{-tv} dv 
\leq e^{-b}/t$ for any $b>0$
and so $\int_0^{b/t} |\sigma(v)| e^{-tv} dv  \gg_\alpha 1/t$
if $b$ is sufficiently large. Taking $N=b/t$ we deduce that if $N$ is
sufficiently large then $\int_0^{N} |\sigma(v)| dv  \gg_\alpha N$,
and so $\lim \sup |\sigma(u)| \gg_\alpha 1$. 

However 
$$
M_0(u,\chi)=\int_0^u \frac{1-\RE \  \chi(v) }{v} dv = \int_1^u \frac{1-\cos \alpha v }{v} dv
= \log u +O_\alpha (1),
$$
so no estimate of the shape $|\sigma(u)| \ll \exp( -\kappa M_0)$ can hold
(with $\kappa>0$), as in Theorem 5.

\subhead 10b.
Corollary 1$'$ is best possible, up to the constant\endsubhead

Assume that Corollary 1$'$ is not best possible, so that if $M=M(u)$
is sufficiently large then $|\sigma(u)| \leq \epsilon Me^{-M}$.

Select $u$ sufficiently large, and choose $\chi(t)=1$ for $t\le 1$, $\chi(t)=i$ for $1< t\le u/2$, and $\chi(t)=0$ for 
$t>u/2$; let $\sigma$ denote the corresponding solution to (1.8).  Next 
we take $\hchi(t)=\chi(t)$ for $t\le u/2$, or $t>u$, 
and for $u/2< t<u$ choose $\hchi(t)$ 
to be a unit vector pointing in the direction of $\overline{\sigma(u-t)}$.  
Let $\hsigma$ denote the corresponding solution to (1.8). 
By definition we have $\hsigma(u-t) =\sigma(u-t)$ in the range $u/2\le t\le u$;
and so $\hchi(t) \hsigma(u-t) = |\sigma(u-t)|$ throughout this range, 
by our choice of $\hchi(t)$. From (1.9) and then this observation we deduce
$$
\hsigma(u) - \sigma(u) = \int_{u/2}^u \frac{\hchi(t)}{t} \hsigma(u-t) dt
=\int_{u/2}^u \frac{|\sigma(u-t)|}{t} dt \geq \frac 1u 
\int_0^{u/2} |\sigma(v)| dv.    \tag{10.1}
$$

Multiplicative functions such as this have been explored 
in some detail in the literature:\ Let $\alpha$ be a complex number with 
Re$(\alpha)<1$, and let $\rho_\alpha$ denote the unique 
continuous solution to $u\rho_\alpha'(u)=-(1-\alpha)\rho_\alpha(u-1)$,
for $u\ge 1$, with the initial condition $\rho_\alpha(u)=1$ for $u\le 1$
(The Dickman-De Bruijn function is the case $\alpha=0$.)
For $\alpha\in [0,1]$, Goldston and McCurley [3] gave
an asymptotic expansion of $\rho_\alpha$.  Their 
proof is in fact valid for all complex $\alpha$ with Re$(\alpha)<1$, and 
shows that when $\alpha$ is not an integer
$$
\rho_\alpha(u) \sim \frac{e^{\gamma(1-\alpha)}}{\Gamma(\alpha)u^{1-\alpha}},
$$
as $u\to \infty$ (Curiously, when $\alpha$ is 
an integer the behaviour of $\rho_\alpha$ is very different;
in fact $\rho_\alpha(u) =1/u^{u+o(u)}$). We have $\sigma(v)=\rho_i(v)$
for $v\leq u/2$, and so in (10.1) we get:
$\hsigma(u) - \sigma(u) = \{ c+o(1)\} \log u /u$ where
$c=e^\gamma/|\Gamma(i)| = 3.414868086\dots$.

Now we note that 
$$
\align
M(u)&=\min_{y\in \Bbb R} \int_0^u \frac{1-\RE \chi(v)e^{-ivy}}{v} dv
=\min_{y\in \Bbb R} \Big(\int_0^1 \frac{1-\cos(vy)}{v} dv +\int_{1}^{u} 
\frac {1-\sin(vy)}{v} dv\Big),\\
&\geq \log u + \min_{y\in \Bbb R} \int_0^y \frac{1-\cos t + \sin t}{t} dt
-\max_{Y\in \Bbb R} \int_0^Y \frac{\sin t}{t} dt
\geq \log u - 1.851937052\dots
\\
\endalign
$$ 
and similarly ${\hat M}(u)\geq \log (u/2) - 1.851937052\dots$.
Let $c'=e^{1.851937052\dots}=6.372150763\dots$.
Therefore $Me^{-M}\leq \{  c'+o(1)\} \log u /u$ and
${\hat M}e^{-{\hat M}} \leq \{  2c'+o(1)\} \log u /u$,
so that 
$|\hsigma(u)|+|\sigma(u)| \geq \{ c+o(1)\} \log u /u \geq 
\{ c/3c'+o(1)\}({\hat M}e^{-{\hat M}} + Me^{-M})$.
Thus either $|\sigma(u)|\geq (5/28) Me^{-M}$ or
$|\hsigma(u)|\geq (5/28) {\hat M}e^{-{\hat M}}$.
This implies the remarks following Corollaries 1 and 1$'$
since $c/3c'> 5/28 > e^\gamma/10$.


\Refs

\ref 
\no 1
\by T. Bonnesen and W. Fenchel 
\book Theorie der konvexen K{\" o}rper
\publ Chelsea, New York
\yr 1949
\endref



\ref
\no 2
\by P.D.T.A. Elliott 
\paper Extrapolating the mean-values of multiplicative functions
\jour Indag. Math
\vol 51
\yr 1989
\pages 409\--420
\endref

\ref
\no 3
\by D.A. Goldston and K.S. McCurley 
\paper Sieving the positive integers by large primes
\jour J. Number Theory
\vol 28
\yr 1988
\pages 94\--115
\endref


\ref
\no 4
\by A. Granville and K. Soundararajan
\paper The spectrum of multiplicative functions
\jour preprint
\endref

\ref
\no 5
\by G. Hal{\' a}sz 
\paper On the distribution of additive and mean-values of multiplicative
functions
\jour Stud. Sci. Math. Hungar
\vol 6 
\yr 1971
\pages 211\--233
\endref

\ref
\no 6
\by G. Hal{\' a}sz
\paper On the distribution of additive arithmetic functions
\jour Acta Arith.
\vol XXVII 
\yr 1975
\pages 143-152
\endref


\ref
\no 7
\by H. Halberstam and H.-E. Richert
\paper On a result of R. R. Hall
\jour J. Number Theory 
\vol 11 
\yr 1979
\pages 76\--89
\endref

\ref
\no 8
\by R.R. Hall
\paper Halving an estimate obtained from Selberg's upper bound method
\jour Acta Arith
\vol 25
\yr 1974
\pages 347\--351
\endref


\ref
\no 9
\by R.R. Hall
\paper A sharp inequality of Hal{\' a}sz type for the mean value of a
multiplicative arithmetic function
\jour Mathematika
\vol 42
\yr 1995
\pages 144\--157
\endref
 
\ref
\no 10
\by R.R. Hall and G. Tenenbaum
\paper Effective mean value estimates for complex multiplicative
functions
\jour Math. Proc. Camb. Phil. Soc.
\vol 110
\yr 1991
\pages 337\--351
\endref

\ref \no 11
\by A. Hildebrand 
\paper A note on Burgess's character sum estimate
\jour C.R. Acad. Sci. Roy. Soc. Canada \vol 8 \yr 1986 \pages 35\--37
\endref

\ref 
\no 12
\by H.L. Montgomery
\paper A note on the mean values of multiplicative functions
\jour Inst.~Mittag-Leffler, (Report \# 17).
\endref

\ref
\no 13
\by H.L. Montgomery and R.C. Vaughan
\paper Hilbert's inequality
\jour J. Lond. Math. Soc. (2)
\vol 8
\yr 1974
\pages 73--82
\endref

\ref 
\no 14
\by H.L. Montgomery and R.C. Vaughan
\paper Mean-values of multiplicative functions
\jour preprint
\endref

\ref 
\no 15
\by A. Wintner 
\book The theory of measure in arithmetical semigroups
\publaddr Baltimore
\yr 1944
\endref



\ref
\no 16
\by E. Wirsing 
\paper Das asymptotische Verhalten von Summen \" uber multiplikative
Funktionen II
\jour Acta Math. Acad. Sci. Hung
\vol 18
\yr 1967
\pages 411\--467
\endref




\endRefs
\enddocument